\documentclass{article}%
\usepackage{latexsym}
\usepackage{amsmath}
\usepackage{graphicx}
\usepackage{amsfonts}
\usepackage{amssymb}%
\newtheorem{theorem}{Theorem}

\newtheorem{conjecture}[theorem]{Conjecture}
\newtheorem{corollary}[theorem]{Corollary}

\newtheorem{definition}[theorem]{Definition}

\newtheorem{lemma}[theorem]{Lemma}
\newtheorem{notation}[theorem]{Notation}

\newtheorem{proposition}[theorem]{Proposition}
\newtheorem{remark}[theorem]{Remark}

\begin{document}

\author{Andrey Todorov\\UC Santa Cruz\\Department of Mathematics\\Santa Cruz, CA 95064\\Institute of Mathematics\\Bulgarian Academy of Sciences\\Sofia, Bulgaria}
\title{Determinants of the Calabi-Yau Metrics on K3 Surfaces, Discriminants, Theta
Lifts and Counting Problems in the A and B Models\\\ \ \ \ \ \ \ \ \ \ \ \ \ \ \ \ \ \ \ \ \ \textit{To Serge Lang with deep
respect}}
\maketitle
\date{}

\begin{abstract}
The Dedekind eta functions plays important role in different branches of
Mathematics and Theoretical Physics. One way to construct Dedekind Eta
function to use the explicit formula (Kroncker limit formula) for the
regularized determinants of the Laplacian of the flat metric acting of (0,1)
forms on elliptic curves. The holomorphic part of the regularized determinant
is the Dedekind eta functions. In this paper we generalized the above approach
to the case of K3 surfaces. We give an explicit formula of the regularized
determinants of the Laplacians of Calabi Yau metrics on K3 Surfaces, following
suggestions by R. Borcherds. The holomorphic part of the regularized
determinants will be the higher dimensional analogue of Dedekind Eta function.

We give explicit formulas for the number of non singular rational curves with
a fixed volume with respect to a Hodge metric in the case of K3 surfaces with
Picard group unimodular even lattice by using the holomorphic part $\exp
\Phi_{3,19}$ of the regularized determinants $\det\Delta_{(0,1)}$.

We gave the combinatorial interpretation of the restriction of the automorphic
form $\exp\Phi_{3,19}$ on the moduli of K3 surfaces with unimodular Picard
lattices in the A and B models. The results obtained in this paper are related
to some results of Bershadsky, Cecotti, Ouguri and Vafa. See \cite{BCOV}.

\end{abstract}
\tableofcontents

\section{Introduction}

The Dedekind eta function%
\[
\eta(\tau)=q^{1/24}%
{\displaystyle\prod\limits_{n=1}^{\infty}}
\left(  1-q^{n}\right)  ,
\]
where $q=e^{2\pi i\tau}$ plays a very important role in different branches of
mathematics. It is closely related to the study of the moduli of elliptic
curves. One way to construct Dedekind Eta function in case of elliptic curves
is to use the explicit formula for the regularized determinants of the
Laplacian of the flat metric acting of (0,1) forms. See \cite{BT05}. The
holomorphic part of the regularized determinant is the Dedekind eta functions.

In this paper we generalized the above approach to the case of K3 surfaces. We
give an explicit formula of the regularized determinants of the Laplacians of
Calabi Yau metrics on the moduli space of Calabi-Yau metrics on the K3
surface, following suggestions by R. Borcherds. The holomorphic part of the
regularized determinants will be the higher dimensional analogue of Dedekind
Eta function.

Next we will review the moduli theory of K3 surfaces. It was A. Weil who
outline the main problems in the study of the moduli of K3 surfaces. See
\cite{W}. The first main result in the study of moduli of K3 surfaces is due
to Shafarevich and Piatetski-Shapiro. See \cite{PS}. They proved the global
Torelli Theorem for polarized algebraic K3 surfaces. Combining the Theorem of
Shafarevich and Piatetski Shapiro with the description of the mapping class
group of K3 surface one obtain that the moduli space $\mathfrak{M}_{k3,n}$ of
polarized algebraic K3 surfaced with a polarization class $e$ such that
$\left\langle e,e\right\rangle =2n>0$ is a Zariski open set in
\[
\Gamma_{K3,n}^{+}\backslash\mathbb{SO}(2,19)/\mathbb{SO}(2)\times
\mathbb{SO}(19),
\]
where $\Gamma_{K3,n}^{+}$ is an index two subgroup in the group of the
automorphisms $\mathcal{O}_{\Lambda_{K3,n}}^{\ast}(\mathbb{Z})$ of the lattice
$H^{2}\left(  \text{M,}\mathbb{Z}\right)  $ which is isomorphic to
\[
\Lambda_{K3,n}:=-2n\mathbb{Z\oplus U}^{2}\mathbf{\oplus}\mathbb{E}%
_{8}(-1)\oplus\mathbb{E}_{8}(-1).
\]
In \cite{To80} it was proved that every point of $\mathbb{SO}%
(3,19)/\mathbb{SO}(2)\times\mathbb{SO}(1,19)$ corresponds to a marked K3
surface. Based on this result in \cite{KT} it was proved that the moduli space
of Ricci flat metrics on K3 surfaces with a fixed volume is isomorphic to
\[
\mathfrak{M}_{KE}:=\Gamma^{+}\backslash\left(  \mathbb{SO}_{0}%
(3,19)/\mathbb{SO}(3)\times\mathbb{SO}(19)-\mathcal{D}_{KE}\right)  ,
\]
where $\Gamma^{+}$ is a subgroup of index 2 in the group of automorphisms of
the group of the automorphisms of the Euclidean lattice $\Lambda
_{K3}=\mathbb{U}^{3}\mathbf{\oplus}\mathbb{E}_{8}(-1)\oplus\mathbb{E}%
_{8}(-1),$ where%
\[
\mathbb{U=}\left(
\begin{array}
[c]{cc}%
0 & 1\\
1 & 0
\end{array}
\right)
\]
is the hyperbolic lattice and $\mathbb{E}_{8}(-1)$ is the standard lattice and
$\mathcal{D}_{KE}$ is the subspace whose points correspond to Ricci flat
metrics on orbifolds. Donaldson proved in \cite{D} that the mapping class
group $\Gamma$ of a K3 surface is\ a subgroup of index 2 in the group of the
automorphisms of the Euclidean lattice $\Lambda_{K3}$.

Mirror Symmetry is based on the observation that there are two different
models A and B in string theory which define one and the same partition
function. The A model is related to the deformation of the K\"{a}hler-Einstein
metrics. The B-model is related to the deformations of complex structures. To
study mirror symmetric on K3 surfaces we need to define a B-field on a K3
surfaces. It is a class of cohomology $\omega_{X}(1,1)\in H^{1,1}%
($X$,\mathbb{C})$ of type $(1,1)$ on a K3 surface X such that
\[%
{\displaystyle\int\limits_{\text{X}}}
\operatorname{Im}\omega\wedge\operatorname{Im}\omega>0.
\]
The moduli space of marked K3 surfaces with a B-field is isomorphic to
$\mathfrak{h}_{4,20}:=\mathbb{SO}_{0}(4,20)/\mathbb{SO}(4)\times
\mathbb{SO}(20).$ Aspinwall and Morrison proved that the moduli space of Super
Conformal Field Theories with supersymmetry (4,4) is described by $\Gamma
_{B}^{+}\backslash\mathfrak{h}_{4,20},$ where $\Gamma_{B}^{+}$ is a subgroup
of index two in $\mathcal{O}(\Lambda_{K3}).$ It is well known that
$\mathfrak{h}_{4,20}$ parametrizes the four-dimensional oriented subspaces in
$\mathbb{R}^{4,20}$ on which the bilinear form is strictly positive. See
\cite{AM}. To a pair (X,$\omega_{X}(1,1))$ of a marked K3 surface with a
B-field $\omega_{X}(1,1)$ we assign a oriented four dimensional subspace
$E_{\text{X,}\omega_{X}(1,1)}$ in
\[
H^{\ast}(\text{X,}\mathbb{Z})\otimes\mathbb{R}=\left(  H^{0}(\text{X,}%
\mathbb{Z})\oplus H^{2}(\text{X,}\mathbb{Z})\oplus H^{4}(\text{X,}%
\mathbb{Z})\right)  \otimes\mathbb{R}%
\]
on which the bilinear form defined by the cup product is positive. We will
assume that $\left(  H^{0}(\text{X,}\mathbb{Z})\oplus H^{4}(\text{X,}%
\mathbb{Z})\right)  =\mathbb{U}_{0}$ and the B-field $\omega_{X}(1,1)$ we will
be identified with
\begin{equation}
(1,-\frac{1}{2}\left(  \omega_{X}(1,1)\wedge\omega_{X}(1,1)\right)  )\in
H^{0}(\text{X,}\mathbb{Z})\oplus H^{2}(\text{X,}\mathbb{Z})\oplus
H^{4}(\text{X,}\mathbb{Z}). \label{ext}%
\end{equation}
From now on we will consider the B-field $\omega_{X}(1,1)$ as defined by
$\left(  \ref{ext}\right)  .$ The four dimensional subspace $E_{\text{X,}%
\omega_{X}(1,1)}$ contains the two dimensional subspace $E_{\omega_{\text{X}}%
}$ spanned by $\operatorname{Re}\omega_{\text{X}}$ and $\operatorname{Im}%
\omega_{\text{X}},$ where $\omega_{\text{X}}$ is the holomorphic two form on X
defined up to a constant and the two dimensional subspace $E_{\omega_{X}%
(1,1)}$ spanned by $\operatorname{Re}\omega_{X}(1,1)$ and $\operatorname{Im}%
\omega_{X}(1,1),$ where $\omega_{X}(1,1)$ is defined by $\left(
\ref{ext}\right)  $. $E_{\omega_{\text{X}}}$ will the orthogonal to
$E_{\omega_{\text{X}}(1,1)}$ in $E_{\text{X,}\omega_{X}(1,1)}.$

Mirror Symmetry is pretty well understood in the case of K3 surfaces. See
\cite{AM}, \cite{dol} and \cite{To93}. The mirror symmetry is exchanging
$E_{\omega_{\text{X}}}$ with $E_{\omega_{\text{X}}(1,1)}.$ Special case of
mirror symmetry of algebraic K3 surfaces was studied in details in \cite{dol}.

In this paper we will consider the moduli space of K3 surfaces with
$B$-fields. We prove the existence of an automorphic form $\exp\left(
\Phi_{4,20}\right)  $ which vanishes on the totally geodesic subspaces that
are orthogonal to $-2$ vectors form following \cite{B97}.

The regularized determinants of the Laplacian of Ricci flat metrics
$\det(\Delta_{KE})$ acting on $(0,1)$ forms will be a function on on the
moduli space of Einstein metric
\[
\mathfrak{M}_{KE}=\mathcal{O}_{\Lambda_{K3}}^{+}\backslash\mathbb{SO}%
_{0}(3,19)/\mathbb{SO}(3)\times\mathbb{SO}_{0}(19).
\]
R. Borcherds suggested that one can compute the determinants of the Laplacians
of Ricci flat metrics explicitly by using the method of the theta lifts. See
\cite{B97}. In this paper we will give an explicit expression of the
regularized determinants of the Laplacians of CY metrics $\det$ as a function
on the moduli space of Einstein metrics $\mathfrak{M}_{KE}.$

The restriction of $\exp\left(  \Phi_{4,20}\right)  $ on the moduli space of
elliptic K3 surfaces with a section
\[
\mathfrak{M}_{ell}:=\Gamma_{ell}\backslash\mathfrak{h}_{2,10}%
\]
vanishes on the discriminant locus
\[
\mathfrak{D}_{ell}\subset\mathfrak{M}_{ell}\subset\Gamma_{B}^{+}%
\backslash\mathfrak{h}_{4,20}%
\]
which is defined by the points orthogonal to $-2$ vectors. The mirror $Y$ of
the elliptic K3 $X$ with the section has Picard group$Pic(Y)=\mathbb{U}%
\oplus\mathbb{E}(-1)\oplus\mathbb{E}(-1).$ $\exp\left(  \Phi_{4,20}\right)  $
restricted on a line $tL$ in the K\"{a}hler cone $K(Y)$ spanned by the
imaginary part $L$ of a Hodge metric, has a Fourier expansion. The Fourier
coefficients \ $a_{n}$ of
\[
\frac{d}{dt}\log\Phi(it)
\]
in front of $\frac{\exp\left(  -int\right)  }{1-\exp\left(  -int\right)  }$
are positive integers and they count the number of rational curves of fix volume.

In the study of moduli of elliptic curves the Dedekind eta function%
\[
\eta(\tau)=q^{1/24}%
{\displaystyle\prod\limits_{n=1}^{\infty}}
\left(  1-q^{n}\right)  ,
\]
where $q=e^{2\pi i\tau}$ plays a very important role. We will point out the
three main properties of $\eta.$

\begin{enumerate}
\item It is well known fact that $\eta^{24}$ is a automorphic form which
vanishes at the cusp. In fact $\eta^{24}$ is the discriminant of the elliptic curve.

\item The Kronecker limit formula gives the explicit relations between the
regularized determinant of the flat metric on the elliptic and $\eta.$

\item The Fourier expansion of $\frac{d}{dt}\log\eta(it)$ are positive
integers which give the number of elliptic curve that that are covering of the
elliptic \ curve $E_{\tau}$ of degree $n.$
\end{enumerate}

By using the results obtained in \cite{BT05} we prove the analogues of the
above properties of the Dedekind eta functions in case of K3 surface for the
restriction of the function $\exp\left(  \Phi_{3,19}\right)  $ on the moduli
space of K3 surfaces with a unimodular Picard lattice. Thus we establish that
$\exp\left(  \Phi_{3,19}\right)  $ is the analogue of the Dedekind eta
function for K3 surfaces.

We also give the combinatorial interpretation of the restriction of the
function $\exp\left(  \Phi_{3,19}\right)  $ on the moduli space of K3 surfaces
with a unimodular Picard lattice in the A-model and for the first time in the
B-model. In the B-model the holomorphic part of the regularized deteriminant
of CY metric counts invariant vanishing calibrated $2-$cycles related to
finite mondromy operators with a given volume when $\operatorname{Im}%
\omega_{Y}$ has integer periods. By invariant vanishing cycles we mean
vanishing invariant cycles under the monodromy that appeared in a families
$\pi:\mathcal{X}\rightarrow D$ over the unit disk such that $\pi^{-1}%
(0)=X_{0}$ has singularities.

We hope that the combinatorial properties of the holomorphic part of the
regularized determinant of CY metric for CY threefolds also holds in the
B-model. It counts the number of invariant calibrated invaraint $3-$cycles of
infinite monodromy.

There are some relations of this paper with the papers \cite{CD} and \cite{CM}.

\subsection{Acknowledgements}

The author want to acknowledge the help and suggestions of Greg Zuckerman. He
proposed to study the behavior of the regularized determinants on the moduli
space of K3 surfaces fibred by elliptic curves with sections.

Special thanks to Jay Jorgenson for his help and comments. I am grateful to
him for introducing me to the ideas of regularized determinants.

I want to thank Jun Li for his interest in this paper and help. Special thanks
to the Center of Mathematical Sciences of Zhe Jiang University and National
Center for Theoretical Sciences, Mathematical Division, National Tsing Hua
University for the financial support during the preparation of the paper.

\section{Symmetric Spaces $\mathfrak{h}_{p,q}:=\mathbb{SO}_{0}%
(p,q)/\mathbb{SO}(p)\times\mathbb{SO}(q)$}

\subsection{Global Flat Coordinates on the Symmetric Space $\mathfrak{h}%
_{p,q}$}

We will need some basic facts about the symmetric space
\[
\mathfrak{h}_{p,q}:=\mathbb{SO}_{0}(p,q)/\mathbb{SO}(p)\times\mathbb{SO}(q).
\]
The following Theorem is standard.

\begin{theorem}
\label{G}Let $\mathbb{R}^{p,q}$ be a $p+q$ dimensional real vector space with
a metric with signature $(p,q).$ There is a one to one correspondence between
points $\tau$ in $\mathfrak{h}_{p,q}$ and all oriented $p-$dimensional
$E_{\tau}$ subspaces in $\mathbb{R}^{p,q}$ on which the intersection form on
$\mathbb{R}^{p,q}$ is strictly positive.
\end{theorem}

\begin{theorem}
\label{G4}Let $\mathbb{R}^{p,q}$ be a $p+q$ dimensional real vector space with
a metric with signature $(p,q).$ Let $E_{\tau_{0}}$ be a $p-$dimensional
subspace in $\mathbb{R}^{p,q}$ such the restriction of the quadratic form on
$E_{\tau_{0}}$ is strictly positive. Let $e_{1},...e_{p}$ be an orthonormal
basis of $E_{\tau_{0}}.$ Let $e_{p+1},...,e_{p+q}$ be orthogonal vectors to
$E_{\tau_{0}}$ such that $\left\langle e_{i},e_{j}\right\rangle =-\delta_{ij}$
for $p+1\leq i,j\leq p+q.$ Let $E_{\tau}$ be any $p-$dimensional subspace in
$\mathbb{R}^{p,q}$ such that the restriction of the quadratic form in
$E_{\tau}$ is strictly positive$.$ Then there exists a basis $\{g_{1}%
(\tau),...,g_{p}(\tau)\}$ in $E_{\tau}$ such that
\begin{equation}
g_{j}(\tau)=e_{j}+\sum_{i=p+1}^{p+q}\tau_{j}^{i}e_{i}. \label{g1}%
\end{equation}

\end{theorem}

\textbf{Proof:} Let
\begin{equation}
f_{i}=\sum_{j=1}^{p}\mu_{i}^{j}e_{j}+\sum_{j=p+1}^{p+q}\lambda_{i}^{j}e_{j}
\label{g4}%
\end{equation}
be an orthonormal basis of $E_{\tau}$ where $1\leq i\leq p$ and $1\leq j\leq
p+q$ $.$ Let $\left(  A_{ij}(\mu)\right)  $ be the $p\times p$ matrix $\left(
\mu_{i}^{j}\right)  $ whose elements $\mu_{i}^{j}$ are defined by the
expression $\left(  \ref{g4}\right)  .$

\begin{lemma}
\label{G41}$\det(A_{ij}(\mu))\neq0.$
\end{lemma}

\textbf{Proof:} Suppose that $\det(A_{ij}(\mu))=0.$ This implies that
$rk(A_{ij}(\mu))<p.$ Thus the rows vectors of the matrix $A_{ij}(\mu)$ are
linearly independent$.$ So we can find constants $a_{i}$ for $i=1,...,q$ such
that at least one of them is non zero and
\begin{equation}%
{\displaystyle\sum\limits_{i=1}^{p}}
a_{i}\left(
{\displaystyle\sum\limits_{j=1}^{p}}
\mu_{i}^{j}e_{j}\right)  =0. \label{g5}%
\end{equation}
Let us consider the vector
\begin{equation}
g(\tau)=%
{\displaystyle\sum\limits_{i=1}^{p}}
a_{i}g_{i}. \label{g6}%
\end{equation}
Combining $\left(  \ref{g5}\right)  $ and $\left(  \ref{g6}\right)  $ we
obtain that
\begin{equation}
g(\tau)=\sum_{j=p+1}^{p+q}\mu_{j}e_{j}. \label{g7}%
\end{equation}
$\left(  \ref{g7}\right)  $ implies that
\begin{equation}
\left\langle g(\tau),g(\tau)\right\rangle =-2\sum_{j=p+1}^{p+q}\left\vert
\mu_{j}\right\vert ^{2}<0. \label{g9}%
\end{equation}
The definition of the vectors $g_{i}(\tau)$ and $\left(  \ref{g6}\right)  $
imply that $g(\tau)$ is a non zero vector in $E_{\tau}.$ Since on $E_{\tau}$
the restriction of the metric is strictly positive we get
\[
\left\langle g(\tau),g(\tau)\right\rangle >0.
\]
Thus we get a contradiction with $\left(  \ref{g9}\right)  .$ Lemma \ref{G41}
is proved. $\blacksquare.$

Theorem \ref{G4} follows directly from Lemma \ref{G41}. $\blacksquare.$

\begin{corollary}
\label{G43}\textit{There is one to one correspondence between the set of all
}$p\times q$\textit{\ matrices (}$\tau_{i}^{j})$ for $1\leq i\leq p$ and
$p+1\leq j\leq p+q$\textit{\ such that the vectors }$g_{i}(\tau)$\textit{\ for
}$i=1,...,p$ \textit{defined by }$\left(  \ref{g1}\right)  $\textit{\ spanned
a }$p-$\textit{dimensional subspace }$E_{\tau}$\textit{\ in }$\mathbb{R}%
^{p,q}$\textit{\ on which the restriction of the quadratic form }$\left\langle
u,v\right\rangle $\textit{\ is strictly positive and the set of points in
}$\mathfrak{h}_{3,19}$\textit{. Thus }$(\tau_{i}^{j})$ define global
coordinates on $\mathfrak{h}_{3,19}.$
\end{corollary}

\subsection{Decomposition of $\mathfrak{h}_{p,q}$}

The following two fact are well known:

\begin{theorem}
\label{Dec1}We have the following decomposition of $\mathfrak{h}%
_{2,p}=\mathbb{R}^{1,p-1}+\sqrt{-1}\mathfrak{h}_{1,p-1}.$
\end{theorem}

\textbf{Proof: }It is a well known fact that $\mathfrak{h}_{1,p-1}$ is one of
the component $V^{+}$ of the cone $V:=\left\{  v\in\mathbb{R}^{1,p-1}%
|\left\langle v,v\right\rangle >0\right\}  .$ Let us consider $\mathbb{R}%
^{2,p}=\mathbb{R}^{1,p-1}\oplus\mathbb{R}^{1,1}.$ Let us consider the map:
\[
\Psi:\mathbb{R}^{1,p-1}+\sqrt{-1}\mathfrak{h}_{1,p-1}\rightarrow
\mathbb{P}\left(  \left(  \mathbb{R}^{1,p-1}\oplus\mathbb{R}^{1,1}\right)
\otimes\mathbb{C}\right)
\]
defined as follows%
\[
\Psi:w=(w_{1},...,w_{p})\rightarrow\left(  w_{1},...,w_{p},-\frac{\left\langle
w,w\right\rangle }{2},1\right)  .
\]
It is easy to check that in $\mathbb{P}\left(  \mathbb{R}^{2,p}\otimes
\mathbb{C}\right)  $ we have
\[
\left\langle \Psi(w),\Psi(w)\right\rangle =0\text{ \& }\left\langle
\Psi(w),\overline{\Psi(w)}\right\rangle >0.
\]
Thus the image of $\mathbb{R}^{1,p-1}+\sqrt{-1}\mathfrak{h}_{1,p-1}$ under the
map $\Psi$ will be $\mathfrak{h}_{2,p},$ since $\mathfrak{h}_{2,p}$ in
$\mathbb{P}\left(  \mathbb{R}^{2,p}\otimes\mathbb{C}\right)  $ is given by one
of the components of the open set in the quadratic $\left\langle
w,w\right\rangle =0$ defined by $\left\langle w,\overline{w}\right\rangle >0.$
It is very easy to prove that $\Psi$ is one to one map. $\blacksquare$

\begin{theorem}
\label{Dec}Suppose that $p\geq3,$ and $q\geq2.$ Then we have the following
decomposition $\mathfrak{h}_{p,q}=\mathfrak{h}_{p-1,q-1}\times\mathbb{R}%
^{p-1,q-1}\times\mathbb{R}_{+},$ where $\mathbb{R}_{+}$ is the set of real
positive numbers$.$
\end{theorem}

\textbf{Proof: }Let us consider in the space $\mathbb{R}^{p,q}$ two vectors
$e_{p+q-1}$ and $e_{p+q}$ such that
\[
\left\langle e_{p+q},e_{p+q}\right\rangle =\left\langle e_{p+q-1}%
,e_{p+q-1}\right\rangle =0\text{ and }\left\langle e_{p+q-1},e_{p+q}%
\right\rangle =1.
\]
Clearly the orthogonal complement to the subspace $\{e_{p+q},e_{p+q}\}$ will
be isometric to $\mathbb{R}^{p-1,q-1}.$ Let us consider a basis $\{e_{1}%
,...,e_{p+1}\}$ of $\mathbb{R}^{p,q},$ where $e_{1},...,e_{p+q-2}$ is a basis
of $\mathbb{R}^{p-1,q-1}.$

There is one to one correspondence between the points $\tau\in\mathfrak{h}%
_{p,q}$ and the oriented $p-$dimensional subspaces $E_{\tau}$ in
$\mathbb{R}^{p,q}$ on which the bilinear form is strictly positive. The
intersection $E_{\tau}\cap\mathbb{R}^{p-1,q-1}$ will be $\left(  p-1\right)
-$dimensional subspace in $\mathbb{R}^{p-1,q-1}$ on which the bilinear form is
strictly positive. Let $f_{1}$ be a vector in $\mathbb{R}^{p,q}$ orthogonal to
$\mathbb{R}^{p-1,q-1}\cap E_{\tau}.$ It is easy to see that the coordinates of
$f_{1}$ can be normalized in such a way that its coordinates in $\mathbb{R}%
^{p,q}$ are such that
\[
f_{1}=(\mu_{1},...,\mu_{p+q-2},1,\lambda),
\]
where $\mu=(\mu_{1},...,\mu_{p+q-2})$ is any vector in $\mathbb{R}^{p-1,q-1}$
and $\lambda>0$ and $\lambda>\left\langle \mu,\mu\right\rangle .$ Thus the
correspondence $E_{\tau}\rightarrow\left(  f_{1},E_{\tau}\cap\mathbb{R}%
^{p-1,q-1}\right)  $ establishes the decomposition $\left(  \ref{F}\right)  .$
$\blacksquare$

\subsection{Definition of the Standard Metric on $\mathfrak{h}_{p,q}$}

Since $\mathfrak{h}_{p,q}\subset Grass(p,p+q)$ then the tangent space
$T_{\tau_{0},\mathfrak{h}_{p,q}}$ at a point $\tau_{0}\in\mathfrak{h}_{p,q}$
can be identified with $Hom\left(  E_{\tau_{0}},E_{\tau_{0}}^{\perp}\right)
.$ Thus any tangent vector $A\in T_{\tau_{0},\mathfrak{h}_{p,q}}$ can be
written in the form%
\begin{equation}
A=%
{\displaystyle\sum\limits_{i=1}^{p}}
{\displaystyle\sum\limits_{j=p+1}^{p+q}}
\tau_{i}^{j}\left(  e_{i}^{\ast}\otimes e_{j}\right)  , \label{M}%
\end{equation}
where $e_{i}$ for $i=1,...,q$ is an orthonormal basis of $E_{\tau_{0}}$ and
$e_{j}$ for $j=p+1,...,p+q$ is an orthonormal basis of \ $E_{\tau_{0}}^{\perp
}.$ Then we define the metric on $T_{\tau_{0},\mathfrak{h}_{p,q}}=Hom\left(
E_{\tau_{0}},E_{\tau_{0}}^{\perp}\right)  $ for $A\in T_{\tau_{0}%
,\mathfrak{h}_{p,q}}=Hom\left(  E_{\tau_{0}},E_{\tau_{0}}^{\perp}\right)  $
given by
\begin{equation}
\left\Vert A^{2}\right\Vert =Tr\left(  A\times A^{t}\right)  =%
{\displaystyle\sum\limits_{i,j}}
\left\vert \tau_{i}^{j}\right\vert ^{2}, \label{met}%
\end{equation}
where $\tau_{j}^{i}$ are defined by $\left(  \ref{M}\right)  .$ We will call
this metric the Bergman metric on $\mathfrak{h}_{3,19}.$

\begin{lemma}
\label{MET}The Bergman metric $ds_{B}^{2}$ is invariant metric on
$\mathfrak{h}_{p,q}.$ It is given in the flat coordinate system $\left(
\tau_{j}^{i}\right)  $ by%
\begin{equation}
ds_{B}^{2}=%
{\displaystyle\sum\limits_{1\leq j\leq3,\text{ }1\leq i\leq19}}
\left(  d\tau_{j}^{i}\right)  ^{2}+O(2). \label{met1}%
\end{equation}

\end{lemma}

\textbf{Proof: }The proof of Lemma \ref{MET} follows directly from the
definition of the Bergman metric. $\blacksquare$

\section{Discriminants in $\mathfrak{H}_{p,q}$}

\subsection{Definition and Basic Properties of the Discriminant}

From now on we will consider the symmetric spaces $\mathfrak{h}_{p,q}$ for
which $p-q\equiv0\operatorname{mod}8.$ In this paper $\Lambda_{p,q}$ will be
unimodular even lattice of signature $(p=q+8k,q).$ We have the following
description all $\Lambda_{p,q}:$

\begin{theorem}
\label{uel}Suppose that $\Lambda_{p,q}$ the unimodular even lattice of
signature $(p,q)$ for $p-q\equiv0\operatorname{mod}8.$ Then
\[
\Lambda_{p,q}\approxeq\underset{p=q+8k}{\underbrace{\mathbb{U}\oplus
...\oplus\mathbb{U}}}\oplus\underset{q}{\underbrace{\mathbb{E}_{8}%
(-1)\oplus...\oplus\mathbb{E}_{8}(-1)}}.
\]

\end{theorem}

\begin{definition}
Define the set $\Delta_{p,q}(e):=\{\delta\in\Lambda_{p,q}|\left\langle
\delta,\delta\right\rangle =-2\}.$ Let us define by $\mathcal{O}_{p,q}$ the
group of the automorphisms of the lattice $\Lambda_{p,q}.$ Let $\mathcal{O}%
_{p,q}^{+}$ be the subgroup of $\mathcal{O}_{p,q}$ which preserve the
orientation of the positive subspaces of dimension $p$ in $\Lambda
_{p,q}\otimes\mathbb{R}.$ Then $\mathcal{O}_{p,q}^{+}$ has index two in
$\mathcal{O}_{p,q}.$
\end{definition}

\begin{definition}
We know that $\mathfrak{h}_{p,q}$ can be realized as an open set in the
Grassmanian $Grass(p,p+q).$ Let us denote by $\mathfrak{h}_{p,q-1}(\delta)$
the set of all $p-$dimensional subspaces in the orthogonal complement of the
vector $\delta$ in $\Lambda_{p,q}\otimes\mathbb{R}.$ We will define the
discriminant locus $\mathfrak{D}_{p,q}$ in $\mathcal{O}_{p,q}^{+}%
\backslash\mathfrak{h}_{p,q}$ as follows:
\[
\mathfrak{D}_{p,q}:=\mathcal{O}_{p,q}^{+}\backslash\left(  \underset{\delta
\in\Delta(e)}{\cup}(\mathfrak{h}_{p,q-1}(\delta))\right)  .
\]

\end{definition}

This definition is motivated by the definition of the discriminant locus in
the moduli of algebraic K3 surfaces.

\subsection{The Irreducibility of the Discriminant}

\begin{theorem}
\label{Bor}The discriminant locus $\mathcal{D}_{p,q}$ is an irreducible real
analytic subspace in $\mathcal{O}_{p,q}^{+}\backslash\mathfrak{h}_{p,q}$,
where $\Lambda_{p,q}$ is an even unimodular lattice.
\end{theorem}

\textbf{Proof:} The proof of Theorem \ref{Bor} will follow if we prove that on
the set of vectors $\Delta_{\Lambda_{p,q}}$ the group $\mathcal{O}%
_{\Lambda_{p,q}}^{+}$ acts transitively. Thus they form one orbit and
therefore the discriminant locus $\mathcal{D}_{p,q}$ in $\mathcal{O}_{p,q}%
^{+}\backslash\mathfrak{h}_{p,q}$ is an irreducible divisor.

The proof that on the set of vectors $\Delta_{\Lambda_{p,q}}$ the group
$\mathcal{O}_{\Lambda_{p,q}}^{+}$ acts transitively will be based on ideas
used in \cite{B96} to prove the irreducibility of the discriminant locus in
the moduli space of Enriques surfaces.

We will proceed by induction on $p$ to prove that the action of $\mathcal{O}%
_{\Lambda_{p,q}}^{+}$ on the set $\Delta_{\Lambda_{p,q}}$ is transitive. For
$p=0$ the Theorem \ref{Bor} is obvious. Suppose that Theorem \ref{Bor} is true
for $p>0.$ We will denote by $L$ the lattice
\[
\underset{p}{\underbrace{\mathbb{U}\oplus...\oplus\mathbb{U}}}\oplus
\underset{q}{\underbrace{\mathbb{E}_{8}(-1)\oplus...\oplus\mathbb{E}_{8}(-1)}}%
\]
and by $M$ the lattice $M=L\oplus\mathbb{U}.$

The plan of the proof is the following. We will denote by $R_{0}$ and $R_{1}$
the set of norm $-2$ vectors of $M$ which have inter product respectively $0$
or $1$ with the vector
\[
e=(\overrightarrow{0},0,1)\in L\oplus\mathbb{U=}M.
\]
Let $\Gamma_{1}$ be the group generated by reflections of elements of the set
$R_{1}$ and $\Gamma_{2}$ be the group generated by reflections of elements of
$R_{0}\cup R_{1}$ and $-id.$ We will show first that any $-2$ vector of $M$ is
conjugate to an element of the set $R_{0}\cup R_{1}.$ Then we will show that
the group $\mathcal{O}_{M}^{+}(\mathbb{Z})$ interchange the sets $R_{0}$ and
$R_{1}.$

\begin{lemma}
\label{Bor3}Any norm $-2$ vector $\delta$ of $M$ is conjugate to an element of
$R_{0}\cup R_{1}$ under the group $\Gamma_{1}.$
\end{lemma}

\textbf{Proof:} The proof of Lemma \ref{Bor3} is based on the following
Propositions \ref{Bor1} and \ref{Bor2}:

\begin{proposition}
\label{Bor1}Suppose that $v\notin M\subset M\otimes\mathbb{Q}.$ Suppose that
$x$ is some real number. Then there exists a vector $\overrightarrow{\mu}\in
M$ such that%
\[
\left\vert \left\langle \overrightarrow{\mu}-\overrightarrow{v}%
,\overrightarrow{\mu}-\overrightarrow{v}\right\rangle -x\right\vert <1.
\]

\end{proposition}

\textbf{Proof:} The proof of Proposition \ref{Bor1} follows the proof of Lemma
\textbf{2.1} given in \cite{B96}. Since $\overrightarrow{v}\notin M\subset
M\otimes\mathbb{Q}$ we can find a primitive isotropic vector $\overrightarrow
{\rho}$ such that $\left\langle \overrightarrow{\rho},\overrightarrow
{v}\right\rangle $ is not an integer. This is because primitive isotropic
vectors span $L.$ As the group $\mathcal{O}_{M}(\mathbb{Z})$ acts transitively
on norm $0$ vectors we can assume that
\[
e=(\overrightarrow{0},0,1)\in\underset{p}{\underbrace{\mathbb{U}%
\oplus...\oplus\mathbb{U}}}\oplus\underset{q}{\underbrace{\mathbb{E}%
_{8}(-1)\oplus...\oplus\mathbb{E}_{8}(-1)}}\oplus\mathbb{U}%
\]
Then $\overrightarrow{v}=(\overrightarrow{\lambda},a,b)$ with $a$ not an
integer. We will find some $\overrightarrow{\mu}$ of the form $\overrightarrow
{\mu}$ $=(\overrightarrow{0},m,n)$ with integers $m$ and $n$ such that for $x$
we have
\[
\left\vert \left\langle \overrightarrow{\mu}-\overrightarrow{v}%
,\overrightarrow{\mu}-\overrightarrow{v}\right\rangle -x\right\vert
=\left\vert \left\langle \overrightarrow{v},\overrightarrow{v}\right\rangle
-2(a-m)(b-n)-x\right\vert <1.
\]
Since $a$ is not an integer we can find some integer $m$ such that $\left\vert
a-m\right\vert <1.$ Whenever we add $1$ to $n,$ the expression $2(a-m)(b-n)$
is changed by a non zero number less than $2,$ so we can choose some integer
$n$ such that $2(a-m)(b-n)$ is at a distance of less then $1$ from any given
number $x-\left\langle \overrightarrow{v},\overrightarrow{v}\right\rangle .$
This proves Proposition \ref{Bor1}. $\blacksquare$

\begin{proposition}
\label{Bor2}Suppose that $R_{1}$ is the set of norm $-2$ vectors of $M$ having
inner product $1$ with
\[
e=(\overrightarrow{0},0,1)\in M=L\oplus\mathbb{U}.
\]
Suppose that $\Gamma_{1}$ is the subgroup of $\mathcal{O}_{M}(\mathbb{Z})$
generated by reflections of vectors of $R_{1}$ and the automorphism $-id.$
Then any vector $r\in M$ is conjugate under $\Gamma_{1}$ to a vector of the
form $(\overrightarrow{v},m,n)\in M$ such that either $m=0$ or $\frac
{\overrightarrow{v}}{m}\in L$ and $m>0.$
\end{proposition}

\textbf{Proof:} We can assume that $\overrightarrow{r}=(\overrightarrow
{v},m,n)$ has the property that $\left\vert \left\langle \overrightarrow
{r},\overrightarrow{e}\right\rangle \right\vert =|m|$ is minimal among all
conjugates of $r$ under $\Gamma_{1},$ where $\overrightarrow{e}%
=(\overrightarrow{0},0,1)$. If $m=0$ then we are done. So we can assume that
$m\neq0,$ and wish to prove that $\frac{v}{m}\in L$ and $m>0.$

Suppose that $\frac{\overrightarrow{v}}{m}\notin L.$ By Proposition \ref{Bor1}
we can find a vector $\overrightarrow{\mu}\in L$ satisfying%
\begin{equation}
\left\vert \left\langle \overrightarrow{\mu}-\frac{\overrightarrow{v}}%
{m},\overrightarrow{\mu}-\frac{\overrightarrow{v}}{m}\right\rangle +\left(
-\frac{2n}{m}-\frac{\left\langle \overrightarrow{v},\overrightarrow
{v}\right\rangle }{m^{2}}\right)  \right\vert <1. \label{b0}%
\end{equation}
Let $\delta=\left(  \overrightarrow{\mu},1,\frac{-\left\langle \overrightarrow
{\mu},\overrightarrow{\mu}\right\rangle -2}{2}\right)  .$ It is easy to see
that $\left\langle \delta,\delta\right\rangle =-2.$ Let
\[
T_{\delta}(\overrightarrow{r})=\overrightarrow{r^{\prime}}=\overrightarrow
{r}+\left\langle \overrightarrow{r},\delta\right\rangle \delta,
\]
i.e. $r^{\prime}$ is the reflection of $r$ with respect to the hyperplane of
$\delta\in R_{1}.$ Direct computations show that%
\[
\left\vert \left\langle \overrightarrow{r^{\prime}},\overrightarrow
{e}\right\rangle \right\vert =\left\vert \left\langle T_{\delta}%
(\overrightarrow{r}),\overrightarrow{e}\right\rangle \right\vert =\left\vert
\left\langle \overrightarrow{r},T_{\delta}(\overrightarrow{e})\right\rangle
\right\vert =\left\vert \left\langle \overrightarrow{r},\overrightarrow
{e}+\left\langle \overrightarrow{e},\delta\right\rangle \delta\right\rangle
\right\vert =
\]%
\begin{equation}
\left\vert m\left(  \left\langle \overrightarrow{\mu}-\frac{\overrightarrow
{v}}{m},\overrightarrow{\mu}-\frac{\overrightarrow{v}}{m}\right\rangle
+\left(  -\frac{2n}{m}-\frac{\left\langle \overrightarrow{v},\overrightarrow
{v}\right\rangle }{m^{2}}\right)  \right)  \right\vert . \label{b2}%
\end{equation}
Combining $\left(  \ref{b0}\right)  $ and $\left(  \ref{b2}\right)  $ we
deduce that%
\begin{equation}
\left\vert \left\langle \overrightarrow{r^{\prime}},\overrightarrow
{e}\right\rangle \right\vert =\left\vert m\left(  \left\langle \overrightarrow
{\mu}-\frac{\overrightarrow{v}}{m},\overrightarrow{\mu}-\frac{\overrightarrow
{v}}{m}\right\rangle +\left(  -\frac{2n}{m}-\frac{\left\langle \overrightarrow
{v},\overrightarrow{v}\right\rangle }{m^{2}}\right)  \right)  \right\vert <m.
\label{b3}%
\end{equation}
We have chosen $\left\vert \left\langle \overrightarrow{r},\delta\right\rangle
\right\vert =m$ to be minimal. So $\left(  \ref{b3}\right)  $ contradicts
$\left\vert \left\langle \overrightarrow{r},\delta\right\rangle \right\vert
=m$. Proposition \ref{Bor2} is proved. $\blacksquare$

\textbf{Proof of Lemma }\ref{Bor3}: Let $\delta=(\overrightarrow{v},m,n).$ By
Proposition \ref{Bor2} we can assume that either $m=0$ or $\frac
{\overrightarrow{v}}{m}\in M.$ If $m=0$ then Lemma \ref{Bor3} is proved.
Suppose that $\frac{\overrightarrow{v}}{m}\in M$ holds. Then $\left\langle
\frac{\overrightarrow{v}}{m},\frac{\overrightarrow{v}}{m}\right\rangle
\in2\mathbb{Z}.$ Thus
\[
\left\langle \delta,\delta\right\rangle =-2=m^{2}\left\langle \frac
{\overrightarrow{v}}{m},\frac{\overrightarrow{v}}{m}\right\rangle +2mn
\]
implies that $-2$ is divisible by $2m.$ From here we conclude that $m=1.$
Lemma \ref{b3} is proved. $\blacksquare$

Let us define the group $\Gamma_{3}$ as the group generated by the
automorphisms $\mathcal{O}_{L}(\mathbb{Z})^{+}$ extended to automorphisms of
$M$ by letting them act trivially on $\mathbb{U},$ the group of automorphisms
taking
\[
(\overrightarrow{v},m,n)\rightarrow(\overrightarrow{v}+2m\overrightarrow
{\lambda},m,n-\left\langle \overrightarrow{v},\overrightarrow{\lambda
}\right\rangle -m\left\langle \overrightarrow{\lambda},\overrightarrow
{\lambda}\right\rangle
\]
for $\lambda\in L,$ and the group of automorphisms given by reflections of
norm $-2$ vectors in $R_{1}.$

\begin{lemma}
\label{Bor4}The group $\Gamma_{3}$ acts transitively on the set of vectors of
norm $-2$ in $M.$
\end{lemma}

\textbf{Proof:} The proof of Lemma \ref{Bor4} is based on the following Propositions:

\begin{proposition}
\label{Bor4a}The group $\mathcal{O}_{L}(\mathbb{Z})^{+}$ acts transitively on
the set of vectors of norm $-2$ in $L.$
\end{proposition}

\textbf{Proof:} Since by definition
\[
L=\underset{p-1}{\underbrace{\mathbb{U}\oplus...\oplus\mathbb{U}}}%
\oplus\underset{q}{\underbrace{\mathbb{E}_{8}(-1)\oplus...\oplus\mathbb{E}%
_{8}(-1)}}\oplus\mathbb{U}%
\]
then Proposition \ref{Bor4a} follows from the induction hypothesis.
$\blacksquare$

\begin{proposition}
\label{Bor4b}Let $\overrightarrow{\lambda}\in L.$ There exists an element
$g_{\overrightarrow{\lambda}}\in Aut\left(  M\right)  $ such that if
$\delta=(\overrightarrow{v},0,k)$ and $\delta^{2}=-2,$ then
$g_{\overrightarrow{\lambda}}(\delta)=(\mu,0,0).$
\end{proposition}

\textbf{Proof:} The conditions $\delta=(\overrightarrow{v},0,k)$ and
$\delta^{2}=-2$ imply that $\left\langle \overrightarrow{v},\overrightarrow
{v}\right\rangle =-2.$ Thus $\overrightarrow{v}$ is a primitive element in
$L.$ Proposition \ref{Bor4a} implies that there exists an element $\sigma
\in\mathcal{O}_{L}(\mathbb{Z})^{+}$ such that
\[
\sigma(\overrightarrow{v})=\overrightarrow{\mu}=f_{1}-\mu_{\overrightarrow
{\mu}}f_{2}\in\mathbb{U\subset}L,
\]
where $\left\langle f_{i},f_{i}\right\rangle =0$ and $\left\langle f_{1}%
,f_{2}\right\rangle =1.$

Next we will construct an element $g_{\overrightarrow{\lambda}}\in Aut\left(
M\right)  $ such that
\[
g_{\overrightarrow{\lambda}}\left(  \left(  \overrightarrow{\mu},0,k\right)
\right)  =(\overrightarrow{v_{1}},0,0).
\]
We know that $\left\langle f_{1},f_{2}\right\rangle =1.$ We can find a
primitive element $\lambda\in L$ such that \textbf{1. }$\left\langle
\overrightarrow{\mu},\overrightarrow{\lambda}\right\rangle \neq0,$ \textbf{2.}
$\overrightarrow{\lambda}\in\mathbb{U}(f_{1},f_{2})^{\perp},$ where
$\mathbb{U}(f_{1},f_{2})$ is the sublattice in $M$ spanned by $f_{1}$ and
$f_{2}$ and \textbf{3. }$\left\langle \overrightarrow{\lambda},\overrightarrow
{v}\right\rangle =k.$ Let us consider the transformation $g_{\overrightarrow
{\lambda}}$ defined as follows:%
\[
g_{\overrightarrow{\lambda}}(f_{2})=f_{2},\text{ }g_{\overrightarrow{\lambda}%
}(\overrightarrow{\alpha})=\overrightarrow{\alpha}-\left\langle \alpha
,\overrightarrow{\lambda}\right\rangle f_{2}\text{ }%
\]
and%
\begin{subequations}
\begin{equation}
g_{\overrightarrow{\lambda}}(f_{1})=f_{1}+\overrightarrow{\lambda}%
-\frac{\left\langle \overrightarrow{\lambda},\overrightarrow{\lambda
}\right\rangle }{2}f_{2}, \label{jb3a}%
\end{equation}
where $\overrightarrow{\alpha}\in M$ is any element. Direct computations show
that $g_{\overrightarrow{\lambda}}$ preserve the scalar product in $M.$ Thus
$g_{\overrightarrow{\lambda}}\in Aut(M).$ The definition of $g_{\lambda}$ and
since we choose $\lambda\in L$ such that $\left\langle \overrightarrow
{\lambda},\overrightarrow{v}\right\rangle =k$ then%
\end{subequations}
\[
g_{\overrightarrow{\lambda}}\left(  \left(  \overrightarrow{\mu},0,k\right)
\right)  =(\overrightarrow{\mu},0,k)-\left\langle \overrightarrow{\mu
},\overrightarrow{\lambda}\right\rangle f_{2}\text{ }=(\overrightarrow
{v},0,0).
\]
Proposition \ref{Bor4b} is proved. $\blacksquare$

\begin{proposition}
\label{Bor4c}Suppose that $\delta\in R_{0}.$ Then there exists an element
$\sigma\in\Gamma_{3}$ such that $\sigma(\delta)\in R_{1}.$
\end{proposition}

\textbf{Proof:} Proposition \ref{Bor4b} implies that without loss of
generality we may assume that $\delta\in L.$ Let
\[
\overrightarrow{\lambda}\in L,\left\langle \delta,\overrightarrow{\lambda
}\right\rangle \neq0\ and\text{ }\left\langle \overrightarrow{\lambda
},\overrightarrow{\lambda}\right\rangle \neq0.
\]
Let us consider%
\[
\delta_{1}=\left(  \overrightarrow{\lambda},1,-\frac{\left\langle
\overrightarrow{\lambda},\overrightarrow{\lambda}\right\rangle +2}{2}\right)
\in L\oplus\mathbb{U}=M.
\]
Clearly $\left\langle \delta_{1},\delta_{1}\right\rangle =-2.$ Then the map
$T_{\delta_{1}}(\delta)=\delta+\left\langle \delta_{1},\delta\right\rangle
\delta_{1}$ is an element of $\Gamma_{3}$ and clearly $T_{\delta_{1}}%
(\delta)\in R_{1}.$ Proposition \ref{Bor4c} is proved. $\blacksquare$

Lemma \ref{Bor3}, Propositions \ref{Bor4a}, \ref{Bor4b} and \ref{Bor4c}
implied Lemma \ref{Bor4}. $\blacksquare$

Lemma \ref{Bor4} implies directly Theorem \ref{Bor}. $\blacksquare$

\section{Moduli of K3 Surfaces}

\subsection{Moduli of Marked, Algebraic and Polarized K3 surfaces}

A K3 surface is a compact, complex two dimensional manifold with the following
properties: \textbf{i. }There exists a non-zero holomorphic two form $\omega$
on X without zeroes. \textbf{ii. }$H^{1}($X,$\mathcal{O}_{\text{X}})=0.$

In \cite{Ast} and \cite{BPV}, the following topological properties are proved.
The surface X is simply connected, and the homology group $H_{2}%
(X,\mathbb{Z})$ is a torsion free abelian group of rank 22. The intersection
form $\left\langle u,v\right\rangle $ on $H_{2}(X,\mathbb{Z})$ has the
properties: \textbf{1.\ }$\left\langle u,u\right\rangle =0$ mod$(2).$
\textbf{2. }$\det\left(  \left\langle e_{i},e_{j}\right\rangle \right)  =-1$.
\textbf{3. }The symmetric form $<$ , $>$ has a signature $(3,19).$

Theorem \textbf{5} on page 54 of \cite{Se} implies that as an Euclidean
lattice $H_{2}(X,\mathbb{Z})$ is isomorphic to the K3 lattice $\Lambda_{K3},$
where $\Lambda_{K3}:=\mathbb{U}^{3}\oplus(-\mathbb{E}_{8})^{2}.$ Every K3
surface is also simply connected.

\begin{definition}
Let $\alpha=\{\alpha_{i}\}$ be a basis of $H_{2}(X,\mathbb{Z})$ with
intersection matrix $\Lambda_{K3}.$ The pair $(X,\alpha)$ is called a marked
K3 surface. Let
\[
l\in H^{1,1}(X,\mathbb{R})\cap H^{2}(X,\mathbb{Z})
\]
be the Poincare dual class of a hyperplane section, i.e. an ample divisor. The
triple $(X,\alpha,l)$ is called a marked, polarized K3 surface. The degree of
the polarization is an integer $2d$ such that $\left\langle l,l\right\rangle
=2d>0.$
\end{definition}

\begin{definition}
The period map $\pi$ for marked K3 surfaces (X,$\alpha)$ is defined by
integrating the holomorphic two form $\omega$ along the basis $\alpha$ of
$H_{2}(X,\mathbb{Z}),$ meaning
\[
\pi(X,\alpha):=(...,\int_{\alpha_{i}}\omega,...)\in\mathbb{P}^{21}.
\]

\end{definition}

The Riemann bilinear relations hold for $\pi(X,\alpha),$ meaning
\begin{equation}
\left\langle \pi(X,\alpha),\pi(X,\alpha)\right\rangle =0\text{ }%
and\ \left\langle \pi(X,\alpha),\overline{\pi(X,\alpha)}\right\rangle >0.
\label{rbr}%
\end{equation}
Choose a primitive vector $l\in\Lambda_{K3}$ such that $\left\langle
l,l\right\rangle =2d>0.$ Let us denote
\[
\Lambda_{K3,l}:=\left\{  v\in\Lambda_{K3}|\left\langle l,v\right\rangle
=0\right\}  .
\]
Then $\pi(X,\alpha,l)\in\mathbb{P}\left(  \Lambda_{K3,l}\otimes\mathbb{C}%
\right)  $ and it satisfies $\left(  \ref{rbr}\right)  .$ The set of points in
$\mathbb{P}\left(  \Lambda_{K3,l}\otimes\mathbb{C}\right)  $ that satisfy
$\left(  \ref{rbr}\right)  $ consists of two components isomorphic to the
symmetric space $\mathfrak{h}_{2,19}.$ In \cite{PS} the following Theorem was proved:

\begin{theorem}
\label{ps}The moduli space $\mathcal{M}_{K3,mpa}^{2d}$ of marked, polarized,
algebraic K3 surfaces of a fixed degree $2d$ exists and it is embedded by the
period map into $\mathfrak{h}_{2,19}$ is an open everywhere dense subset. Let
\[
\Gamma_{K3,2d}=\{\phi\in Aut^{+}(\Lambda_{K3})|\left\langle \phi
(u),\phi(u)\right\rangle =\left\langle u,u\right\rangle \text{ }and\text{
}\phi(l)=l\},
\]
where $l$ is a primitive vector such that $\left\langle l,l\right\rangle
=2d>0.$ Then the moduli space $\mathcal{M}_{K3,pa}^{2d}$ of polarized,
algebraic K3 surfaces of a fixed degree 2d is isomorphic to a Zariski open set
in the quasi-projective variety $\Gamma_{K3,2d}$%
$\backslash$%
$\mathfrak{h}_{2,19}.$
\end{theorem}

By pseudo-polarized algebraic K3 surface we understand a pair (X,$l$) where
$l$ corresponds to either ample divisor or pseudo ample divisor, which means
that for any effective divisor $D$ in X, we have $\left\langle
D,l\right\rangle \geq0.$ Mayer proved the linear system $|3l|$ defines a map:
\[
\phi_{|3l|}:X\rightarrow X_{1}\subset\mathbb{P}^{m}%
\]
such that: \textbf{i. }$X_{1}$ has singularities only double rational points.
\textbf{ii. }$\phi_{|3l|}$ is a holomorphic birational map. Let us denote by
\ $\mathcal{M}_{K3,ppa}^{2d}$ the moduli space of pseudo-polarized algebraic
K3 surfaces of degree 2d. From the results proved in \cite{D}, \cite{Ku},
\cite{To80} and \cite{PS} the following Theorem follows:

\begin{theorem}
\label{pskt}The moduli space of $\mathcal{M}_{K3,ppa}^{2d}$ is isomorphic to
the locally symmetric space $\Gamma_{K3,2d}\backslash\mathfrak{h}_{2,19}.$
\end{theorem}

\subsection{Moduli Space of K3 Surfaces with a B-Field}

\begin{definition}
\label{b-field}Let $X$ be a K3 surface. Let $\omega_{X}(1,1)\in H^{1,1}\left(
X\text{,}\mathbb{C}\right)  $ such that
\[%
{\displaystyle\int\limits_{X}}
\operatorname{Im}\omega_{X}(1,1)\wedge\operatorname{Im}\omega_{X}(1,1)>0.
\]
Then $\omega_{X}(1,1)$ will be called a B-field on X.
\end{definition}

\begin{theorem}
\label{mb-field}Let $\left(  X,\omega_{X}\left(  1,1\right)  ,\gamma
_{1},...,\gamma_{22}\right)  $ be a marked K3 surface with a B-field. Then the
moduli space $\mathfrak{M}_{m,B}$ of marked K3 surfaces with a B-field is
isomorphic to $\mathbb{SO}_{0}\left(  4,20\right)  /\mathbb{SO}\left(
4\right)  \times\mathbb{SO}\left(  20\right)  .$
\end{theorem}

\textbf{Proof: }See \cite{AM} or \cite{To93}. $\blacksquare$

\subsection{Discriminant of Pseudo-Polarized K3 Surfaces}

The complement of $\mathcal{M}_{K3,mpa}^{2d}$ in $\mathfrak{h}_{2,19}$ can be
described as follow. Given a polarization class $e\in\Lambda_{K3},$ set
$T_{e}$ to be the orthogonal complement to $e$ in $\Lambda_{K3},$ i.e. $T_{e}$
is the transcendental lattice. Then we have the realization of $\mathfrak{h}%
_{2,19}$ as one of the components of
\[
\mathfrak{h}_{2,19}\approxeq\{u\in\mathbb{P}(T_{e}\otimes\mathbb{C}%
)|\left\langle u,u\right\rangle =0and\left\langle u,\overline{u}\right\rangle
>0\}.
\]
For each $\delta\in\Delta(e),$ define the hyperplane
\[
H(\delta)=\{u\in\mathbb{P}(T_{e}\otimes\mathbb{C})|\left\langle u,\delta
\right\rangle =0\}.
\]
Let $\mathcal{H}_{K3,2d}=\underset{\delta\in\Delta(e)}{\cup}(H(\delta
)\cap\mathfrak{h}_{2,19}).$ Let us define the discriminant $\mathcal{D}%
_{K3}^{2d}:=\Gamma_{K3,2d}\backslash\mathcal{H}_{K3,2d}.$ Results from
\cite{Ma}, \cite{PS}, \cite{To80} and \cite{Ku} imply that $\mathcal{D}%
_{K3}^{2d}$ is the complement of the moduli space of algebraic polarized K3
surfaces $\mathcal{M}_{K3,pa}^{2d}$ in the locally symmetric space
$\Gamma_{K3,2d}\backslash\mathfrak{h}_{K3,2d},$ i.e.
\[
\mathcal{D}_{K3}^{2d}=(\Gamma_{K3,2d}\backslash\mathfrak{h}_{K3,2d}%
)\ -\mathcal{M}_{K3,pa}^{2d}.
\]

\section{Mirror Symmetry}

\subsection{Mirror Symmetry for K3 Surfaces}

Let $(X,\alpha,\omega_{X}(1,1))$ be a marked K3 surface with a B-field
$\omega_{X}(1,1)$. To define the mirror of $(X,\alpha,\omega_{X}(1,1))$ we
need to fix an unimodular hyperbolic lattice $\mathbb{U}$ in $H_{2}%
(X,\mathbb{Z})$ with generators $\left\{  \gamma_{0},\gamma_{1}\right\}  $
such that for the holomorphic two form $\omega_{\text{X}}$ we have
\[%
{\displaystyle\int\limits_{\gamma_{0}}}
\omega_{X}\neq0\text{ }and%
{\displaystyle\int\limits_{\gamma_{1}}}
\omega_{X}\neq0.
\]
Thus we can normalize $\omega_{X}$ in the following manner
\begin{equation}%
{\displaystyle\int\limits_{\gamma_{0}}}
\omega_{X}=1\text{ and }%
{\displaystyle\int\limits_{\gamma_{1}}}
\omega_{X}\neq0. \label{ms0}%
\end{equation}
From now on we will consider the set $(X,\alpha,\omega_{X}(1,1),\mathbb{U}%
,\mathbb{\omega}_{X}),$ where $\alpha$ is a marking, $\mathbb{U}$ is a fixed
sublattice in $H^{2}(X,\mathbb{Z})$ such that the holomorphic two form
satisfies $\left(  \ref{ms0}\right)  .$ Let $\mathbb{U}^{\perp}$ be the
orthogonal complement of $\mathbb{U}$ in $H^{2}(X,\mathbb{Z}).~$Let us denote
by $\mathbb{U}_{0}$ the unimodular hyperbolic sublattice $H^{0}(X,\mathbb{Z}%
)\oplus H^{4}(X,\mathbb{Z})$ in the cohomology ring $H^{\ast}(X,\mathbb{Z}).$
We will assign to the $B$-field $\omega_{X}(1,1)$ the vector%
\[
\hat{\omega}_{X}:=\left(  \omega_{X}(1,1),1,-\frac{\omega_{X}(1,1)\wedge
\overline{\omega_{X}(1,1)}}{2}\right)
\]
in $H^{2}(X,\mathbb{Z})\oplus\mathbb{U}_{0}=H^{\ast}(X,\mathbb{Z}).$

We will need the following Theorem:

\begin{theorem}
\label{mir}Let $(X,\alpha,\omega_{X}(1,1),\mathbb{U},\mathbb{\omega}_{X}),$
where $\alpha$ is a marking, $\mathbb{U}$ is a fixed sublattice in
$H^{2}(X,\mathbb{Z})$ such that the holomorphic two form satisfies $\left(
\ref{ms0}\right)  .$ Then there exists a marked K3 surface $(Y,\alpha)$ with a
$B$-field $\omega_{Y}(1,1)$ such that \textit{if we identify }$H^{2}%
(Y,\mathbb{Z})$\textit{ with} $\mathbb{U}^{\perp}\oplus\mathbb{U}_{0},$
\textit{then the class of the cohomology }$[\omega_{Y}]$\textit{ of the K3
surfaces }$Y$ is such that $[\omega_{Y}]=\omega_{X}(1,1)\mathit{\in}\left(
\mathbb{U}^{\perp}\oplus\mathbb{U}_{0}\right)  \otimes\mathbb{C}$ and\textit{
}$\omega_{Y}(1,1)=[\omega_{X}]\mathit{\in}\left(  H^{2}(Y,\mathbb{Z}%
)\oplus\mathbb{U}\right)  \otimes\mathbb{C}.$
\end{theorem}

\textbf{Proof: }Let us consider $\omega_{X}(1,1)\in\left(  \mathbb{U}^{\perp
}\oplus\mathbb{U}_{0}\right)  \otimes\mathbb{C=}\Lambda_{K3}\otimes
\mathbb{C}.$ Then direct computations show that we have $\left\langle
\hat{\omega}_{X},\hat{\omega}_{X}\right\rangle =0$ and $\left\langle
\hat{\omega}_{X},\overline{\hat{\omega}_{X}}\right\rangle >0.$ From the
epimorphism of the period map for K3 surfaces proved in \cite{To80} it follows
that there exists a marked K3 surface $\,(Y,\alpha)$ with a holomorphic two
form $\omega_{Y}$ such that the class of cohomology $[\omega_{Y}]$ is the same
as the class of cohomology of $\omega_{X}(1,1).$ Next we will prove that the
class of cohomology $\omega_{X}\in H^{1,1}(Y,\mathbb{C})$ satisfies
\[%
{\displaystyle\int\limits_{Y}}
\operatorname{Im}\omega_{X}\wedge\operatorname{Im}\omega_{X}=\left\langle
\operatorname{Im}\omega_{X},\operatorname{Im}\omega_{X}\right\rangle >0.
\]
Indeed on $X$ we have
\begin{equation}
\left\langle \omega_{X},\omega_{X}(1,1)\right\rangle =\left\langle \omega
_{X},\overline{\omega_{X}(1,1)}\right\rangle =0 \label{mir0}%
\end{equation}
since $\omega_{X}(1,1)$ is a form of type $(1,1)$ and $\omega_{X}$ is a form
of type $(2,0).$ On the other hand the form $\omega_{X}(1,1)$ with respect to
the new complex structure $Y$ on X it is a form of type $(2,0\dot{)}.$ So
$\left(  \ref{mir0}\right)  $ means that on $Y$ $\omega_{X}$ is a form of type
$(1,1\dot{)}.$On the other hand we have%
\begin{equation}%
{\displaystyle\int\limits_{X}}
\omega_{X}\wedge\overline{\omega_{X}}=2%
{\displaystyle\int\limits_{X}}
\operatorname{Im}\omega_{X}\wedge\operatorname{Im}\omega_{X}=2\left\langle
\operatorname{Im}\omega_{X},\operatorname{Im}\omega_{X}\right\rangle >0.
\label{mir1}%
\end{equation}
Thus $\left(  \ref{mir1}\right)  $ proves that $\omega_{X}$ is a B-field on
$Y.$ Theorem \ref{mir} is proved. $\blacksquare$

Now we are ready to define the mirror symmetry:

\begin{definition}
\label{MS}We will define the marked surface $(Y,\alpha,\omega_{Y}%
(1,1),\mathbb{U},\omega_{Y})$ constructed in Theorem \ref{mir} will be the
mirror of $(X,\alpha,\omega_{X}(1,1),\mathbb{U}_{0},\omega_{X}).$
\end{definition}

\subsection{Mirror Symmetry and Algebraic K3 Surfaces}

Let us consider the Neron-Severi group
\[
M=Pic(X):=H^{2}(X,\mathbb{Z})\cap H^{1,1}(X,\mathbb{R}).
\]
We can characterize in another way $NS(X).$ It is the dual group in
$H^{2}(X,\mathbb{Z})$ of the kernel of the functional:
\[
\left(  \omega_{X}\right)  :H_{2}(X,\mathbb{Z})\rightarrow\mathbb{C}%
\]
defined by $\gamma\rightarrow%
{\displaystyle\int\limits_{\gamma}}
\omega_{X}.$ We define the transcendental classes of homologies $T(X)\subset$
$H^{2}(X,\mathbb{Z})$ on $X$ as follows: $T(X):=\ker\left(  \omega_{X}\right)
^{\perp}.$

\begin{definition}
We will say that pairs $(X,$M$)$ M$-$marked K3 surface if M is the Picard
lattice of some algebraic K3 surface together with a primitive imbedding of M
into $H^{2}(X,\mathbb{Z})$.
\end{definition}

The following Theorem was proved in \cite{To93} or \cite{dol}

\begin{theorem}
\label{mmod}The moduli space $\mathfrak{M}_{\text{M}}$ of marked pairs
$(X,$M$)$ exists and $\mathfrak{M}_{\text{M}}\approxeq\Gamma_{\text{M}%
}\left\backslash \mathfrak{h}_{2,20-\rho}\right.  ,$ where $\rho=rk$M and
$\Gamma_{\text{M}}=\left\{  \phi\in Aut\Lambda_{K3}\left\vert \phi\right.
_{\text{M}}=id\right\}  .$
\end{theorem}

Suppose that we consider $M$ such that $\mathbb{U}$ can be embedded into
$M^{\perp}.$ According to a Theorem of Nikulin this is always possible if
$rkM=\rho\geq9.$ The construction of mirror symmetry for $M$ marked K3
surfaces $(X,\alpha,M,\mathbb{U},\omega_{X}),$ where $\mathbb{U}\subset
M^{\perp}$ was described in \cite{To93} and \cite{dol} as follows; Let
$(X,\alpha,M,\mathbb{U},\omega_{X})$ be an algebraic polarized K3 surface.
Then Theorem \ref{mir} implies the following Corollary:

\begin{corollary}
\label{mira}Let $(X,\alpha,M,\mathbb{U},\omega_{X})$ be M-marked K3 surface
such that $\mathbb{U}\subset T_{X}$ and the $B$-field $\omega_{X}(1,1)$
satisfies $\omega_{X}(1,1)|_{\mathbb{U}^{\perp}\subset T_{X}}=0.$ Then the
mirror $(Y,M_{1},\mathbb{U},\omega_{Y}(1,1),\omega_{Y})$ satisfies the
following conditions: \textbf{i. }$Pic(Y)=M_{1}=\mathbb{U}^{\perp}\subset
T_{X}.$ \textbf{ii. }$T_{Y}=M\oplus\mathbb{U}\approxeq Pic(X)\oplus
\mathbb{U}.$
\end{corollary}

\textbf{Proof: }Corollary \ref{mira} follows directly from Theorem \ref{mir}.
$\blacksquare$

\begin{remark}
Some interesting examples and applications of Corollary \ref{mira} were
discussed in \cite{dol}.
\end{remark}

\subsection{The Mirror Map for Marked M-K3}

Part of the mirror conjecture states that the

\begin{definition}
Let X be a K3 surface. We will define the K\"{a}hler cone of $K$(X) of X as
follows:%
\[
K\text{(X)}:=\left\{  \omega\in H^{1,1}\left(  \text{X,}\mathbb{R}\right)
|\omega=\operatorname{Im}g,\text{ and }g\text{ is a K\"{a}hler metric on
X}\right\}  .
\]

\end{definition}

We will need the characterization of the K\"{a}hler cone that is given bellow.
Denote by
\[
\Delta(X):=\left\{  \delta\in NS(X)|\left\langle \delta,\delta\right\rangle
=-2\right\}  .
\]
We will need the following Lemma from \cite{PS}:

\begin{lemma}
Let $\delta\in\Delta(X).$ Then $\delta$ or $-\delta$ can be realized as an
effective curve on X.
\end{lemma}

We will denote by
\[
\Delta^{+}(X):=\left\{  \delta\in\Delta(X)|\delta\text{ can be realized as an
effective cure}\right\}  .
\]
Let us denote by $V:=\left\{  v\in H^{1,1}\left(  X,\mathbb{R}\right)
|\left\langle v,v\right\rangle >0\right\}  .$ Since the restriction of the
bilinear form on $H^{1,1}\left(  X,\mathbb{R}\right)  $ has a signature
$(1,19),$ then $V$ will consists of two components. Let us denote by $V^{+}$
the component of $V$ which contains a K\"{a}hler class.

Each $\delta\in$ $\Delta^{+}(\Delta)$ generates a reflection $s_{\delta}$ of
$V^{+},$ where $s_{\delta}(v)=v+\left\langle v,\delta\right\rangle \delta.$
Let us denote by $\Gamma(\Delta)$ the subgroup of $\mathcal{O}_{\Lambda_{K3}%
}^{+}$ generated by $s_{\delta}.$ In \cite{To80} the following Theorem was proved:

\begin{theorem}
\label{to89}The K\"{a}hler cone $K(X)$ coincides with the fundamental domain
of the group $\Gamma(\Delta)$ in $V^{+}$ which contains a K\"{a}hler class.
\end{theorem}

\textbf{Proof: }See \cite{To80}. $\blacksquare$

\begin{remark}
\label{mirrm}According to Theorem \ref{mmod} $\mathfrak{M}_{K3,M}=\Gamma
_{M}\backslash\mathfrak{h}_{2,20-\rho}$ is the moduli space of $M$-marked K3
surfaces. Suppose that $\mathbb{U}\subset T_{X}$ is fixed and $M_{1}\subset
T_{X}$ is the orthogonal complement of $\mathbb{U}$ in $M$. Let $(Y,M_{1})$ be
some $M_{1}$ marked K3 surface defined by the primitive embedding
$M_{1}\subset T_{X}\subset\Lambda_{K3}.$ Let $\mathfrak{h}_{M_{1}}%
=M_{1}\otimes\mathbb{R}+iK(Y),$ where $K(Y)$ is the K\"{a}hler cone of $Y.$
Then according to Theorem \ref{Dec1} $\mathfrak{M}_{K3,M}\approxeq\Gamma
_{M}\backslash\mathfrak{h}_{M_{1}}.$ Thus we have a complex analytic covering
map:%
\[
\psi_{M}:\mathfrak{h}_{M_{1}}\rightarrow\mathfrak{M}_{K3,M}=\Gamma
_{M}\backslash\mathfrak{h}_{M_{1}}.
\]
The map $\psi_{M}^{-1}$ which is multivalued is called the mirror map. It
identifies in the case described in this Remark the moduli space of $M$-marked
K3 surfaces with the complexified K\"{a}hler cone of the its mirror.
\end{remark}

\section{Automorphic Forms on $\Gamma\backslash\mathfrak{h}_{p,q},$ Theta
Lifts and Regularized Determinants of CY metrics on K3 Surfaces}

\subsection{General Facts about Regularized Determinants}

\begin{definition}
Let M be a compact $C^{\infty}$ manifold. Let $g$ be a Rimannian metric on M.
Let
\[
\Delta_{g,q}=d\circ d^{\ast}+d^{\ast}\circ d
\]
be the Laplacian associated with the metric $g$ acting on the space of
$C^{\infty}$ infinity $q-$forms $C^{\infty}\left(  \text{M,}\Omega_{\text{M}%
}^{q}\right)  .$ It is a well known fact that the spectrum of $\Delta_{g,q}$
is non negative, i.e. $0\leq\lambda_{1}\leq...\lambda_{k}\leq$ ...and
\begin{equation}
\underset{k\rightarrow\infty}{\lim}\frac{\lambda_{k}}{k^{\frac{2}{n}}}=c>0,
\label{lim}%
\end{equation}
where $n=\dim_{\mathbb{R}}$M. We define the zeta function $\zeta_{q}(s)$ of
$\Delta_{g,q}$ as follows:%
\[
\zeta_{q}(s)=%
{\displaystyle\sum\limits_{\lambda_{k}>0}}
\lambda_{k}^{-s}.
\]
Then $\left(  \ref{lim}\right)  $ implies that $\zeta_{q}(s)$ is a well
defined function for $s\in\mathbb{C}$, where $\operatorname{Re}s$ large
enough. One can prove that $\zeta_{q}(s)$ has a meromorphic continuation in
$\mathbb{C}$ and $\zeta_{q}(0)$ is well defined. Then we define the
regularized determinant of $\Delta_{g,q}$ as follows: $\det\Delta_{g,q}%
=\exp\left(  -\zeta_{q}^{\prime}(0)\right)  .$
\end{definition}

In \cite{BT05} the following Theorem was proved:

\begin{theorem}
\label{BT}Let M be a CY manifold with a polarization class $L\in H^{2}%
($M$,\mathbb{Z})\cap H^{1,1}($M$,\mathbb{R}).$ Let $\det\Delta_{(0,1)}$ be the
regularized determinant of the Laplacian corresponding to the Calabi Yau
metric corresponding to the polarization class $L$ and acting on the space of
$(0,1)$ forms. Then $dd^{c}\log\det\Delta_{(0,1)}=-\operatorname{Im}W.P..$
\end{theorem}

\subsection{Special Automorphic Form of Weight -2 on $\Gamma\backslash
\mathfrak{h}_{p,q}$}

In this paper the group $\Gamma$ will be the group of automorphisms of
$\Lambda_{K3}$ which preserve the spinor norm, i.e. $\Gamma=\mathcal{O}%
_{\Lambda_{K3}}^{+}(\mathbb{Z})$ is a subgroup of index $2$ in the group of
automorphism $\mathcal{O}_{\Lambda_{K3}}(\mathbb{Z})$ of the lattice
$\Lambda_{K3}.$ Donaldson proved in \cite{D} that the mapping class group of a
K3 surface is isomorphic to $\Gamma.$

We will define the one cocycle $\mu(\gamma,\tau)$ of the group $\Gamma$ with
coefficients the non singular $3\times3$ matrices with coefficients functions
on $\mathfrak{h}_{3,19}.$ Let an element $\gamma\in\Gamma$ be represented by a
$(22\times22)$ matrix $(\gamma_{k,l}).$ We proved that any point $\tau
\in\mathfrak{h}_{3,19}$ can be represented by the $3\times22$ matrix
$\tau=(E_{3},\tau_{ij}),$ where $E_{3}$ is the identity $3\times3$
matrix$.$The action of $\gamma=(\gamma_{k,l})\in\Gamma$ on $\mathfrak{h}%
_{3,19}$ is described as follow:
\begin{equation}
\gamma(\tau)=\left(  E_{3},\tau_{i,j})\times(\gamma_{k,l})\right)
=(\mu(\gamma,\tau),\sigma_{\gamma,ij}(\tau)), \label{d2}%
\end{equation}
where $\mu(\gamma,\tau)$ is $3\times3$ matrix defined by the first three
columns of the matrix $\left(  \ref{d2}\right)  $ and $\sigma_{\gamma,ij}%
(\tau)$ is some $3$ $\times19$ matrix. Theorem \ref{G4} implies that the
$3\times3$ matrix $\mu(\gamma,\tau)$ has rank $3,$ i.e. $\det(\mu(\gamma
,\tau))\neq0.$ It is easy to see that we have:%
\[
\mu(\gamma_{1}\gamma_{2},\tau)=\mu(\gamma_{1},\tau)\times\mu(\gamma_{2}%
,\gamma_{1}(\tau)).
\]

\begin{definition}
\label{G6}Let $\Phi(\tau)$ be a function on $\mathfrak{h}_{3,19}$ such that it
satisfies the following functional equation:
\[
\Phi(\tau\gamma)=(\det\mu(\gamma,\tau))^{k}\Phi(\tau).
\]
Then we will call $\Phi(\tau)$ an automorphic form of weight $k.$
\end{definition}

\begin{definition}
\label{G8} Let us recall that according to Theorem \ref{G4} to each point
$\tau=(\tau_{j}^{i})\in\mathfrak{h}_{3,19},$ $1\leq j\leq3$ and $1\leq
i\leq19$ we assigned the row vectors $g_{i}$ of the matrix ($E_{3},\tau
_{i}^{j}).$ We will define the function $g(\tau)$ on $\mathfrak{h}_{3,19}$ as
follows
\[
g(\tau):=\det(\left\langle g_{i}\left(  \tau\right)  ,g_{j}\left(
\tau\right)  \right\rangle \dot{)}.
\]

\end{definition}

\begin{theorem}
\label{G9}The function $g(\tau)$ defined in Definition \ref{G8} is an
automorphic form of weight $-2.$
\end{theorem}

\textbf{Proof}: We need to compute
\[
g((\gamma(\tau))=\det(\left\langle \left(  \mu(\gamma,\tau)\right)
^{-1}\times\gamma_{i}(\tau)),\left(  \mu(\gamma,\tau)\right)  ^{-1}%
\times\gamma_{j}(\tau))\right\rangle )=?,
\]
where $\gamma_{i}(\tau)$ is the $i^{th}$ row of the $(3\times22)$ matrix
$(\tau_{j}^{i})\times\gamma.$ Theorem \ref{G4} and the expression of the
matrix $\mu(\gamma,\tau)$ given by $\left(  \ref{d2}\right)  $ imply
\[
g((\gamma(\tau))=\det(\left\langle \left(  \mu(\gamma,\tau)\right)
^{-1}\times\gamma_{i}(\tau)),\left(  \mu(\gamma,\tau)\right)  ^{-1}%
\times\gamma_{j}(\tau)))\right\rangle =
\]%
\[
\left(  \det(\mu(\gamma,\tau))\right)  ^{-2}\det(\left\langle g_{i}%
(\tau),g_{j}(\tau)\right\rangle \dot{)}=\det(\mu(\gamma,\tau))^{-2}\times
g(\tau).
\]
Thus we get $g(\gamma(\tau))=\det(\mu(\gamma,\tau))^{-2}g(\tau).$ So Theorem
\ref{G9} is proved. $\blacksquare$

\subsection{Theta Lifst and Automorphic Form with a Zero Set Supported by the
Discriminant Locus on $\Gamma\backslash\mathfrak{h}_{3,19}$}

Suppose that $\Lambda_{p,q}$ is a unomodular even lattice. We will define the
Siegel kernel $\Theta_{\Lambda_{p,q}}(\tau)$ as follows:%
\[
\Theta_{\Lambda_{p,q}}(\tau):=%
{\displaystyle\sum\limits_{\lambda\in\Lambda_{p,q}}}
\exp\left(  2\pi\sqrt{-1}\left\langle \Pr\text{ }_{E_{\tau}}\lambda,\Pr\text{
}_{E_{\tau}}\lambda\right\rangle \rho-\left\langle \Pr\text{ }_{E_{\tau}%
^{\bot}}\lambda,\Pr\text{ }_{E_{\tau}^{\bot}}\lambda\right\rangle
\overline{\rho}\right)  ,
\]
where $E_{\tau}$ is a $p-$dimensional real vector subspace in $\Lambda
_{p,q}\otimes\mathbb{R}$ on which the quadratic form is positive definite,
$E_{\tau}^{\bot}$ is the $q-$dimensional vector subspace in $\Lambda
_{p,q}\otimes\mathbb{R}$ orthogonal to $E_{\tau},$ $\Pr$ $_{E_{\tau}}\lambda$
is the orthogonal projection of $\lambda$ on $E_{\tau},$ $\Pr$ $_{E_{\tau
}^{\bot}}\lambda$ is the orthogonal projection of $\lambda$ on $E_{\tau}%
^{\bot},$and $\rho=x+iy,$ $y>0.$

The following result follows directly from the results proved in \cite{B97}.

\begin{theorem}
\label{Bor5}Let $\Lambda_{p,q}$ be an even unimodular lattice of signature
$(p,q).$ Then there exists a non zero automorphic form $\exp\left(
\Phi_{\Lambda_{p,q}}(\tau)\right)  $ such that the zero set of $\exp\left(
\Phi_{\Lambda_{p,q}}(\tau)\right)  $ coincide with the discriminant
\[
\mathcal{D}_{\Lambda_{p,q}}\mathcal{\subset O}_{\Lambda_{p,q}}^{+}%
\backslash\mathfrak{h}_{p,q}.
\]
Moreover let $\Lambda_{p_{1},q_{1}}$ be an even unimodular sublattice in
$\Lambda_{p,q}$ such that $p-q=p_{1}-q_{1}.$ Then
\begin{equation}
\exp\left(  \Phi_{\Lambda_{p,q}}(\tau)\right)  \left\vert _{\mathcal{O}%
_{\Lambda_{p_{1},q_{1}}}^{+}(\mathbb{Z})\backslash\mathfrak{h}_{p_{1},q_{1}}%
}\right.  =\exp\left(  \Phi_{\Lambda_{p_{1},q_{1}}}(\tau)\right)  .
\label{res0}%
\end{equation}

\end{theorem}

\textbf{Proof: }Let us consider the regularized integral as described in
\cite{B97} or in \cite{K}
\[
\Phi_{\Lambda_{p,q}}(\tau)=%
{\displaystyle\int\limits_{\mathcal{H}}}
\Theta_{\Lambda_{p,q}}(\tau)y^{\frac{q}{2}}\frac{E(\rho)}{\Delta(\rho)}%
\frac{d\rho\wedge\overline{d\rho}}{y^{2}},
\]
where $\mathcal{H}$ is the fundamental domain of the group $\mathbb{PSL}%
_{2}(\mathbb{Z}),$ $\rho=x+iy,$ $y>0$ and $\frac{E(\rho)}{\Delta(\rho)}$ is a
meromorphic automorphic form of weight $q-p$ with a pole of order one at
$\infty.$ Thus we have
\begin{equation}
\frac{E(\rho)}{\Delta(\rho)}=\frac{1}{\exp\left(  2\pi\sqrt{-1}\rho\right)
}+a_{0}+a_{1}\exp\left(  2\pi\sqrt{-1}\rho\right)  +... \label{res1}%
\end{equation}
It was proved in \cite{B97} and in \cite{K} that $\left(  \ref{res1}\right)  $
implies that $\exp\left(  \Phi_{\Lambda_{p,q}}(\tau)\right)  $ will vanish on
the discriminant of $\mathcal{O}^{\ast}(\Lambda_{p,q})\backslash
\mathfrak{h}_{p,q}.$

The relation $\left(  \ref{res0}\right)  $ follows from the condition
$p-q=p_{1}-q_{1}=8k$ and the definition of $\Phi_{\Lambda_{p,q}}(\tau).$
Theorem \ref{Bor5} is proved. $\blacksquare$

We will consider the case of K3 surfaces. We know that $\Lambda_{K3}%
=\Lambda_{3,19}.$ We will study the relations between the non zero automorphic
form $\exp\left(  \Phi_{\Lambda_{K3}}(\tau)\right)  $ and the regularized determinants.

\begin{theorem}
\label{a11}$\Delta_{B}\Phi_{\Lambda_{K3}}(\tau,\sigma)=0,$ where $\Delta_{B}$
is the Laplacian of the Bergman metric on $\Lambda_{K3}=\Lambda_{3,19}.$
\end{theorem}

\textbf{Proof:} Any choice of an embedding of the hyperbolic lattice
$\mathbb{U\subset}\Lambda_{K3}$ defines a totally geodesic subspace
$\mathfrak{h}_{2,18}$ into $\mathfrak{h}_{3,19}.$ This follows from Theorem
\ref{Dec}. According to the construction of the automorphic form $\exp\left(
\Phi_{\Lambda_{2,18}}(\tau)\right)  $ given in \cite{B97} it follows that
$\Phi_{\Lambda_{ell}}$ is a holomorphic function on $\mathfrak{h}_{2,18}.$
Thus we have $\Delta_{B}\Phi_{\Lambda_{ell}}=0.$ All the embeddings
$\mathfrak{h}_{2,18}\subset\mathfrak{h}_{3,19}$ corresponding to primitive
embeddings $\mathbb{U\subset}\Lambda_{K3}$ form an everywhere dense subset of
totally geodesic subamnifolds in $\mathfrak{h}_{3,19}.$Since $\mathfrak{h}%
_{2,18}$ is a totally geodesic subspace in $\mathfrak{h}_{3,19}$ we get that
\[
\Delta_{B}\left(  \Phi_{\Lambda_{K3}}\left\vert _{\mathfrak{h}_{2,18}}\right.
\right)  =\Delta_{B}\Phi_{\Lambda_{ell}}.
\]
Thus the restriction of the Bergman Laplacian applied to on $\Phi
_{\Lambda_{ell}}$ is zero on an everywhere dense subset in $\mathfrak{h}%
_{3,19}.$ Thus the continuous function $\Delta_{B}\Phi_{\Lambda_{K3}}$ is zero
on everywhere dense subset in $\mathfrak{h}_{3,19}.$ From here we deduce that
$\Delta_{B}\Phi_{\Lambda_{K3}}=0.$ Theorem \ref{a11} is proved. $\blacksquare$

\subsection{The Analogue of the Kronecker Limit Formula for the Regularized
Determinants on K3 Surfaces}

\begin{theorem}
\label{har}The function $\log\frac{\det\Delta_{KE}}{\det\left(  \left\langle
g_{i}(\tau),g_{j}(\tau)\right\rangle \right)  }$ is a harmonic function on the
moduli space $\mathfrak{M}_{\text{KE}}$ of Einstein metrics of the K3 surface
with respect to the Laplacian corresponding to the Bergman metric.
\end{theorem}

\textbf{Proof: }The proof of Theorem \ref{har} is based on the following Lemmas:

\begin{lemma}
\label{har1}Let $\tau_{0}\in\mathfrak{h}_{3,19}.$ Then there exists
$L\in\Lambda_{K3}\otimes\mathbb{R}$ and totally geodesic subspace
$\mathfrak{h}_{2,19}$ passing through $\tau_{0}\in\mathfrak{h}_{3,19}$ such
its points correspond to polarized marked K3 surfaces with class of
polarization $L$.
\end{lemma}

\textbf{Proof: }We know that each point $\tau=(\tau_{j}^{i})\in\mathfrak{h}%
_{3,19}$ corresponds to a three dimensional subspace $E_{\tau}\subset
\Lambda_{K3}\otimes\mathbb{R}$ on which the cup product is strictly positive.
Let $L\in E_{\tau}$ be a non zero vector. Then $\left\langle L,L\right\rangle
>0.$ Let us consider the following set:%
\[
\mathfrak{h}_{L}:=\left\{  \text{ }E\subset\Lambda_{K3}\otimes\mathbb{R}%
\left\vert L\in E,\text{ }\dim_{\mathbb{C}}E=3\text{ \& }\left\langle \text{
},\text{ }\right\rangle \left\vert _{E}\right.  >0\right.  \right\}  \text{ .}%
\]
It is easy to see that there is one to one correspondence between the two
dimensional oriented positive subspaces in the orthogonal complement $L^{\bot
}=\mathbb{R}^{2,19}$ and $\mathfrak{h}_{L}$. Thus we get that
\[
\mathfrak{h}_{L}=\mathfrak{h}_{2,19}=\mathbb{SO}_{0}(2,19)/\mathbb{SO}%
(2)\times\mathbb{SO}(19).
\]
Lemma \ref{har1} is proved. $\blacksquare$

Let us choose an orthonormal basis $e_{1},$ $e_{2}$ and $e_{3}=L$ of the three
dimensional subspace $E_{\tau_{0}}\in\mathfrak{h}_{L}.$ Lemma \ref{har1} and
Corollary \ref{G43} imply that the three dimensional subspaces $E_{\tau}$ that
correspond $\tau\in\mathfrak{h}_{L}=\mathfrak{h}_{2,19}\subset\mathfrak{h}%
_{3,19}$ are spanned by the orthonormal vectors:%
\begin{equation}
g_{1}(\tau)=e_{1}+%
{\displaystyle\sum\limits_{i=1}^{19}}
\tau_{1}^{i}e_{i},\text{ }g_{2}(\tau)=e_{2}+%
{\displaystyle\sum\limits_{i=1}^{19}}
{}\tau_{2}^{i}e_{i}\text{ and }g_{3}(\tau)=L=e_{3}. \label{g8}%
\end{equation}

\begin{lemma}
\label{har2}The subspace given by the equations $\tau_{3}^{i}=0$ for
$i=1,...,19,$ where $\tau_{j}^{i}$ are coordinates defined by $\left(
\ref{g8}\right)  $ is the totally geodesic subspace in $\mathfrak{h}%
_{L}=\mathfrak{h}_{2,19}$ in $\mathfrak{h}_{3,19}$.
\end{lemma}

\textbf{Proof: }The proof follows directly from $\left(  \ref{g8}\right)  .$
$\blacksquare$

We know that $\mathfrak{h}_{2,19}$ is a complex manifold of dimension $19$.
The complex coordinates on $\mathfrak{h}_{2,19}$ are defined as follows:
$\rho^{i}=\tau_{1}^{i}+\sqrt{-1}\tau_{2}^{i},$ $1\leq i\leq19.$ From the
epimorphism of the period map we know that $\tau_{0}$ corresponds to a K3
surface $X_{\tau_{0}}$ and the class of cohomology of the complex two form
$\omega_{\tau_{0}}:=e_{1}+\sqrt{-1}e_{2}\in\Lambda_{K3}\otimes\mathbb{C}$ can
be identified with the class of cohomology of the holomorphic two form
$\omega_{\tau_{0}}(2,0)$ on a marked K3 surface $X_{\tau_{0}}$ such that the
vector $e_{3}=L$ can be identified with the class of cohomology of the
imaginary part of a K\"{a}hler metric on $X_{\tau_{0}}.$ The subspace in
$\Lambda_{K3}\otimes\mathbb{R}$ spanned by $e_{4},...,e_{22}$ can be
identified with the primitive class of cohomology of type $(1,1),$ i.e. with
$H_{0}^{1,1}(X_{\tau_{0}},\mathbb{R})=E_{\tau_{0}}^{\perp}.$ See \cite{To80}.

\begin{definition}
We will define the Weil-Petersson metric on the totally geodesic subspace
$\mathfrak{h}_{2,19}$ as the restriction of the metric on $\mathfrak{h}%
_{3,19}$ defined by $\left(  \ref{met}\right)  .$
\end{definition}

\begin{lemma}
\label{har4}Let $\tau_{0}\in\mathfrak{h}_{3,19}.$ Let $\mathfrak{h}%
_{L}=\mathfrak{h}_{2,19}$ be the totally geodesic subspace passing through
$\tau_{0}\in\mathfrak{h}_{3,19}$ and defined by the $L\in E_{\tau_{0}}$ as in
Lemma \ref{har1}$.$ Let $g_{i}(\tau)$ be vectors defined by $\left(
\ref{g8}\right)  .$ Then the function
\[
\log\det\left(  \left\langle g_{i}(\tau),g_{j}(\tau)\right\rangle \right)
|_{\mathfrak{h}_{2,19}}%
\]
is a potential of the Weil-Petersson metric on $\mathfrak{h}_{2,19}.$
\end{lemma}

\textbf{Proof: }The $2\times2$ matrix $\left(  \left\langle g_{i}(\tau
),g_{j}(\tau)\right\rangle \right)  |_{\mathfrak{h}_{2,19}}$ is symmetric.
Since
\[
\left\langle g_{i}(0),g_{j}(0)\right\rangle =\delta_{ij}%
\]
it can be represented as follows:
\[
\left(  \left\langle g_{i}(\tau),g_{j}(\tau)\right\rangle \right)
|_{\mathfrak{h}_{2,19}}=I_{2}+\left(  h_{ij}(\tau)\right)  .
\]
Then we have:%
\begin{equation}
\log\det\left(  \left\langle g_{i}(\tau),g_{j}(\tau)\right\rangle \right)
\left\vert _{\mathfrak{h}_{2,19}}\right.  =%
{\displaystyle\sum\limits_{i=1}^{2}}
\log(1+\lambda_{i}(\tau)), \label{g10}%
\end{equation}
where $\lambda_{i}\left(  \tau\right)  $ are the eigen values of the matrix
$\left(  h_{ij}(\tau)\right)  .$ Thus we get%
\begin{equation}%
{\displaystyle\sum\limits_{i=1}^{2}}
\lambda_{i}\left(  \tau\right)  =h_{11}(\tau)+h_{22}(\tau). \label{g11}%
\end{equation}
From the definition of the matrix $\left(  \left\langle g_{i}(\tau),g_{j}%
(\tau)\right\rangle \right)  |_{\mathfrak{h}_{2,19}}$ we get that
\begin{equation}
h_{11}=%
{\displaystyle\sum\limits_{i=4}^{22}}
\left(  \tau_{1}^{i}\right)  ^{2}\text{ and }h_{22}=%
{\displaystyle\sum\limits_{i=4}^{22}}
\left(  \tau_{2}^{i}\right)  ^{2}. \label{g12}%
\end{equation}
Combining $\left(  \ref{g10}\right)  ,$ $\left(  \ref{g11}\right)  $ and
$\left(  \ref{g12}\right)  $ we get that
\begin{equation}
\log\det\left(  \left\langle g_{i}(\tau),g_{j}(\tau)\right\rangle \right)
\left\vert _{\mathfrak{h}_{2,19}}\right.  =\frac{1}{4}%
{\displaystyle\sum\limits_{i=4}^{22}}
|\rho^{i}|^{2}+O(3). \label{g14}%
\end{equation}
Thus we get from $\left(  \ref{g14}\right)  $ that
\begin{equation}
dd^{c}\log\det\left(  \left\langle g_{i}(\tau),g_{j}(\tau)\right\rangle
\right)  \left\vert _{\mathfrak{h}_{2,19}}\right.  =\frac{\sqrt{-1}}{2}%
{\displaystyle\sum\limits_{i=4}^{22}}
\partial\rho^{i}\wedge\overline{\partial\rho^{i}}+O(2). \label{g15}%
\end{equation}
From $\left(  \ref{g15}\right)  $ we conclude the proof of Lemma \ref{har4}.
$\blacksquare$

\begin{lemma}
\label{har5}Let $\Delta_{B}$ be the Laplacian of the Bergman metric on
$\mathfrak{h}_{3,19}.$ Then the restriction of the function
\[
\log\frac{\det\left(  \Delta_{KE}\right)  }{\det\left(  \left\langle
g_{i}(\tau),g_{j}(\tau)\right\rangle \right)  }%
\]
on each totally geodesic subspace $\mathfrak{h}_{2,19}\subset\mathfrak{h}%
_{3,19}$ is harmonic function with respect of the Laplacian of the
Weil-Petersson metric.
\end{lemma}

\textbf{Proof: }Combining Theorem \ref{BT} with Lemma \ref{har4} we deduce
Lemma \ref{har5}. $\blacksquare$

It is an obvious fact that the set of three dimensional positive subspaces in
$\Lambda_{K3}\otimes\mathbb{R}$ which contain a vector in $\Lambda_{K3}%
\otimes\mathbb{Q}$ form an everywhere dense subset in $\mathfrak{h}_{3,19}.$
From here it follows that we can find an everywhere dense subset of totally
geodesic subsets $\mathfrak{h}_{2,19}$ in $\mathfrak{h}_{3,19}$ on which the
continuous function%
\[
\Delta_{B}\left(  \log\det\Delta_{KE}-\log\det\left(  \left\langle g_{i}%
(\tau),g_{j}(\tau)\right\rangle \right)  \right)
\]
is zero. Therefore it is zero on $\mathfrak{h}_{3,19}.$ Theorem \ref{har} is
proved. $\blacksquare$

\begin{theorem}
\label{a1}The following formula holds for the regularized determinant of the
Laplacian of the Einstein metrics
\[
\det(\Delta_{KE})(\tau)=\det\left(  \left\langle g_{i}(\tau),g_{j}%
(\tau)\right\rangle \right)  \times\left\vert \exp\left(  \Phi_{\Lambda_{K3}%
}(\tau)\right)  \right\vert ^{2}.
\]

\end{theorem}

\textbf{Proof:} According to Theorem \ref{har} the function
\[
\log\det\Delta_{KE}-\log\det\left(  \left\langle g_{i}(\tau),g_{j}%
(\tau)\right\rangle \right)
\]
is a harmonic function with respect to the Laplacian of the Bergman metric on
$\mathfrak{h}_{3,19}.$ Let us consider the function:
\[
\frac{\det\Delta_{KE}}{\det\left(  \left\langle g_{i}(\tau),g_{j}%
(\tau)\right\rangle \right)  }=\phi\left(  \tau\right)
\]
on $\mathfrak{h}_{3,19}.$ According to Theorem \ref{G9} the function
$\det\left(  \left\langle g_{i}(\tau),g_{j}(\tau)\right\rangle \right)  $ is
an automorphic form of weight $-2.$ Therefore the function $\phi$ is an
automorphic function of weight $2.$ In \cite{JT0} we proved that $\det
\Delta_{KE}$ is a bounded non negative function. Therefore the only zeroes of
$\det\Delta_{KE}$ can be located on the discriminant locus $\mathfrak{D}%
_{KE}.$ We know that $\left\vert \exp\left(  \Phi_{\Lambda_{K3}}(\tau)\right)
\right\vert $ is an automorphic function with a zero set on the discriminant
locus $\mathfrak{D}_{KE}.$ Since $\mathfrak{D}_{KE}$ is an irreducible in
$\mathfrak{M}_{KE},$ by taking suitable powers of $\phi$ and $\left\vert
\exp\left(  \Phi_{\Lambda_{K3}}(\tau)\right)  \right\vert ,$ we may assume
that the function
\[
\frac{\left\vert \exp\left(  \Phi_{\Lambda_{K3}}(\tau)\right)  \right\vert
^{n}}{\phi^{m}}=\psi
\]
is a non zero function such $\Delta_{B}\log\psi=0.$ Thus we get a harmonic non
zero function on $\mathfrak{M}_{KE}.$

\begin{lemma}
\label{a1a} $\psi|_{\mathfrak{M}_{ell}}=const.$
\end{lemma}

\textbf{Proof: }Since
\[
dd^{c}\left(  \log\frac{\det(\Delta_{KE})(\tau)}{\det\left(  \left\langle
g_{i}(\tau),g_{j}(\tau)\right\rangle \right)  }\left\vert \mathfrak{M}%
_{ell}\right.  \right)  =0
\]
we can conclude that%
\[
\frac{\det(\Delta_{KE})(\tau)}{\det\left(  \left\langle g_{i}(\tau),g_{j}%
(\tau)\right\rangle \right)  }\left\vert \mathfrak{M}_{ell}\right.  =|\eta|,
\]
where $\eta$ is a holomorphic automorphic form defined up to a character
$\chi\in\Gamma_{ell}\left/  [\Gamma_{ell},\Gamma_{ell}]\right.  $ and with a
zero set $\mathfrak{D}_{ell}.$ Since $\mathfrak{D}_{ell}$ is an irreducible
divisor, we can conclude that $\eta=\exp\left(  \Phi_{\Lambda_{ell}}%
(\tau)\right)  .$ Thus since%
\[
\exp\left(  \Phi_{\Lambda_{K3}}(\tau)\right)  \left\vert _{\mathfrak{M}_{ell}%
}\right.  =\exp\left(  \Phi_{\Lambda_{ell}}(\tau)\right)  ,
\]
we get that $\psi\left\vert _{\mathfrak{M}_{ell}}\right.  =conct.$ Since any
two $\mathfrak{M}_{ell,1}$ and $\mathfrak{M}_{ell,2}$ intersect. So the
continuous function $\psi$ is a constant on an everywhere dense subset in
$\mathfrak{M}_{KE}$. Thus $\psi$ is a constant. Lemma \ref{a1a} is proved.
$\blacksquare$

Lemma \ref{a1a} imply Theorem \ref{a1}. $\blacksquare$

\section{Harvey-Moore-Borcherds Products and Counting Problems in A and B
Models}

\subsection{Counting Problems on K3}

\begin{theorem}
\label{to05}Let $X$ be an algebraic K3 surface such that Pic(X) is an
unimodular lattice. Then we have either $NS(X)=\mathbb{U\oplus E}_{8}(-1)$ or
$NS(X)=\mathbb{U\oplus E}_{8}(-1)\mathbb{\oplus E}_{8}(-1).$ Let $l\in NS(X)$
be the polarization class. Let us consider the components $V_{Enr}^{+}$ and
$V_{ell}^{+}$ of the positive cones in $\left(  \mathbb{U\oplus E}%
_{8}(-1)\right)  \otimes\mathbb{R}$ and in $\left(  \mathbb{U\oplus E}%
_{8}(-1)\mathbb{\oplus E}_{8}(-1\right)  \otimes\mathbb{R}$ which contain the
polarization vector $l.$ Let us consider the discriminant automorphic forms
$\exp\left(  \Phi_{enr}\left(  \tau\right)  \right)  $ and $\exp\left(
\Phi_{ell}\left(  \tau\right)  \right)  $ on $\left(  \mathbb{U\oplus E}%
_{8}(-1)\right)  \otimes\mathbb{R\oplus}\sqrt{-1}V_{Enr}^{+}$ and on $\left(
\left(  \mathbb{U\oplus E}_{8}(-1)\right)  \otimes\mathbb{R}\right)
\mathbb{\oplus}\sqrt{-1}V_{ell}^{+}.$ Then the restriction of the functions
$\exp\left(  \Phi_{enr}\left(  \tau\right)  \right)  $ and $\exp\left(
\Phi_{ell}\left(  \tau\right)  \right)  $ on the lines $\sqrt{-1}lt$ are
periodic. The Fourier expansions%
\[
\frac{d}{dt}\left(  \Phi_{Enr}\left(  \sqrt{-1}lt\right)  \right)  =-%
{\displaystyle\sum\limits_{n}}
a_{n}\frac{e^{-nt}}{1-e^{-nt}}%
\]
and
\begin{equation}
\text{ }\frac{d}{dt}\left(  \Phi_{ell}\left(  \sqrt{-1}lt\right)  \right)  =-%
{\displaystyle\sum\limits_{n}}
b_{n}\frac{e^{-nt}}{1-e^{-nt}}\text{ } \label{CP}%
\end{equation}
have integer coefficients $a_{n}$ and $b_{n}.$ $a_{n}$ and $b_{n\text{ }}$are
equal to the number of non singular rational curves of degree $n$ on a K3
surface $X$ with $NS(X)=\mathbb{U\oplus E}_{8}(-1)$ or $NS(X)=\mathbb{U\oplus
E}_{8}(-1)\mathbb{\oplus E}_{8}(-1).$
\end{theorem}

\textbf{Proof: }Let us fix a bases $\left\{  \gamma_{i}\right\}  $ and
$\left\{  \varepsilon_{j}\right\}  $ of
\[
\mathbb{U\oplus E}_{8}(-1)\text{ and }\mathbb{U\oplus E}_{8}(-1)\mathbb{\oplus
E}_{8}(-1)
\]
respectively. Then we fix the flat coordinates $\left\{  \tau^{1}%
,...,\tau^{10}\right\}  $ and $\left\{  \tau^{1},...,\tau^{18}\right\}  $ in
the symmetric spaces $\mathfrak{h}_{2,10}$ and $\mathfrak{h}_{2,18}$
respectfully represented as tube domains since we have
\[
\mathfrak{h}_{2,10}=\left(  \mathbb{U\oplus E}_{8}(-1)\right)  \otimes
\mathbb{R}+iV^{+}\subset\left(  \mathbb{U\oplus E}_{8}(-1)\right)
\otimes\mathbb{C}%
\]
and
\[
\mathfrak{h}_{2,18}=\left(  \mathbb{U\oplus E}_{8}(-1)\mathbb{\oplus E}%
_{8}(-1)\right)  \otimes\mathbb{R}+iV^{+}\subset\left(  \mathbb{U\oplus E}%
_{8}(-1)\mathbb{\oplus E}_{8}(-1)\right)  \otimes\mathbb{C},
\]
where $V^{+}$ is one of the components of the positive cone in $\left(
\mathbb{U\oplus E}_{8}(-1)\right)  \otimes\mathbb{R}$ or $\left(
\mathbb{U\oplus E}_{8}(-1)\mathbb{\oplus E}_{8}(-1)\right)  \otimes
\mathbb{R}.$

We will denote by $\left\langle \delta,\tau\right\rangle $ the following
expressions:%
\[
\left\langle \delta,\tau\right\rangle =%
{\displaystyle\sum\limits_{i=1}^{10}}
\left\langle \delta,\gamma_{i}\right\rangle \tau^{i}\text{ and }\left\langle
\delta,\tau\right\rangle =%
{\displaystyle\sum\limits_{i=1}^{18}}
\left\langle \delta,\varepsilon_{i}\right\rangle \tau^{i}.
\]
Then Harvey-Moore-Borcherds product formula states that there exist
automorphic forms on $\Gamma_{2,10}\backslash\mathfrak{h}_{2,10}$ or on
$\Gamma\backslash\mathfrak{h}_{2,18}$ which can be represented for some large
$\operatorname{Im}\tau^{i}$ as the following products.
\[
\exp\left(  \Phi_{Enr}\left(  \tau\right)  \right)  =\exp(2\pi i\left\langle
\tau,w\right\rangle
{\displaystyle\prod\limits_{\delta\in\Delta_{Enr}^{+}}}
\left(  1-\exp\left(  2\pi i%
{\displaystyle\sum\limits_{i=1}^{10}}
\left\langle \delta,\gamma_{i}\right\rangle \tau^{i}\right)  \right)
\]
and
\begin{equation}
\exp\left(  \Phi_{ell}\left(  \tau\right)  \right)  =\exp(2\pi i\left\langle
\tau,w\right\rangle
{\displaystyle\prod\limits_{\delta\in\Delta_{Enr}^{+}}}
\left(  1-\exp\left(  2\pi i%
{\displaystyle\sum\limits_{i=1}^{18}}
\left\langle \delta,\varepsilon_{i}\right\rangle \tau^{i}\right)  \right)  .
\label{cp0}%
\end{equation}
It was proved that $\exp\left(  \Phi_{Enr}\left(  \tau\right)  \right)  $ and
$\exp\left(  \Phi_{ell}\left(  \tau\right)  \right)  $ have an analytic
continuation in $\mathfrak{h}_{2,10}$ and $\mathfrak{h}_{2,18}$ and the zeroes
remain the same. Substituting
\[%
{\displaystyle\sum\limits_{i=1}^{10}}
\gamma_{i}\tau^{i}=ilt\text{ and }%
{\displaystyle\sum\limits_{i=1}^{18}}
\varepsilon_{i}\tau^{i}=ilt\text{ }%
\]
in $\left(  \ref{cp0}\right)  $ we get
\[
\exp\left(  \Phi_{Enr}\left(  \tau\right)  \right)  =\exp(2\pi i\left\langle
\tau,w\right\rangle
{\displaystyle\prod\limits_{\delta\in\Delta_{Enr}^{+}}}
\left(  1-\exp\left(  -2\pi\left\langle \delta,l\right\rangle t\right)
\right)
\]
and
\begin{equation}
\exp\left(  \Phi_{ell}\left(  \tau\right)  \right)  =\exp(2\pi i\left\langle
\tau,w\right\rangle
{\displaystyle\prod\limits_{\delta\in\Delta_{Enr}^{+}}}
\left(  1-\exp\left(  -2\pi\left\langle \delta,l\right\rangle t\right)
\right)  . \label{cp1}%
\end{equation}
Let us split the irreducible non singular on disjoint finite sets $A_{n},$
where $A_{n}=\left\{  \delta\in\Delta^{+}|\left\langle \delta,l\right\rangle
=n\right\}  .$ Suppose that $\#A_{n}=a_{n}$ in the case of $\Lambda_{Enr\text{
}}$ and $\#A_{n}=b_{n}$ in the case $\Lambda_{ell}.$ We can rewrite $\left(
\ref{cp1}\right)  $ as follows
\[
\exp\left(  \Phi_{Enr}\left(  \tau\right)  \right)  =\exp(2\pi i\left\langle
\tau,w\right\rangle
{\displaystyle\prod\limits_{\delta\in\Delta_{Enr}^{+}}}
\left(  1-\exp\left(  -2\pi\left\langle \delta,l\right\rangle t\right)
\right)  =
\]%
\[
\exp(2\pi i\left\langle \tau,w\right\rangle
{\displaystyle\prod\limits_{n=1}}
\left(
{\displaystyle\prod\limits_{\delta\in A_{n}}}
\left(  1-\exp\left(  -2\pi nt\right)  \right)  \right)  =
\]%
\begin{equation}
\exp(2\pi i\left\langle \tau,w\right\rangle
{\displaystyle\prod\limits_{n=1}}
\left(  \left(  1-\exp\left(  -2\pi nt\right)  \right)  ^{a_{n}}\right)  .
\label{cp2}%
\end{equation}
In the same way we will get that%
\[
\exp\left(  \Phi_{ell}\left(  \tau\right)  \right)  =\exp(2\pi i\left\langle
\tau,w\right\rangle
{\displaystyle\prod\limits_{\delta\in\Delta_{Enr}^{+}}}
\left(  1-\exp\left(  -2\pi\left\langle \delta,l\right\rangle t\right)
\right)  =
\]%
\begin{equation}
\exp(2\pi i\left\langle \tau,w\right\rangle
{\displaystyle\prod\limits_{n=1}}
\left(  \left(  1-\exp\left(  -2\pi nt\right)  \right)  ^{b_{n}}\right)  .
\label{cp3}%
\end{equation}
From $\left(  \ref{cp2}\right)  $ and $\left(  \ref{cp3}\right)  $ we derive
$\left(  \ref{CP}\right)  $ and thus Theorem \ref{to05}. $\blacksquare$

\subsection{A and B Models}

\begin{remark}
In the $A$ model the automorphic function $\exp\left(  \Phi_{3,19}\left(
\tau\right)  \right)  $ which is the holomorphic part of the regularized
determinant when restricted on the line $\mathbb{R}l$ in the K\"{a}hler cone
of a K3 surface with $Pic(X)$ unimodular lattice counts rational curves with a
given volume according to Theorem \ref{to05}.
\end{remark}

We will consider the $B$ model of $M-$marked K3 surfaces where $M$ is an
unimodular lattice and $M=Pic(Y).$ The moduli space of $Pic(Y)-$marked K3
surfaces $\mathfrak{M}_{Pic(Y)},$ where $Pic(Y)$ is a unimodular lattice can
be represented as a tube domain $\mathbb{R}^{k}+iV^{+}$ modulo action of an
arithmetic group $\Gamma_{Pic(Y)}$. Now we will study the combinatorial
properties of the restriction of the automorphic function $\exp\left(
\Phi_{4,20}\left(  \tau\right)  \right)  $ on $\mathfrak{M}_{Pic(Y)}.$

\begin{definition}
\label{calcycles}Let $Y$ is a K3 surface. Let $g$ is a Calabi-Yau metric on
$Y.$ Let $\gamma\in H_{2}\left(  Y,\mathbb{Z}\right)  .$ We will call $\gamma$
a calibrated cycle if the restriction of $\alpha\operatorname{Re}\omega
_{Y}+\beta\operatorname{Im}\omega_{Y}$ on $\gamma$ is the volume form of the
restriction of the CY metric on $\gamma.$
\end{definition}

\begin{theorem}
\label{Bmodel}Suppose that Y is a K3 surface such that $\operatorname{Im}%
\omega_{Y}\in H^{2}(Y,\mathbb{Z})\cap H^{1,1}(Y,\mathbb{R}).$ Let $g$ be a CY
metric on Y. Then any $\delta\in T(Y):=Pic(Y)^{\bot}\subset H^{2}%
(Y,\mathbb{Z})$ such that $\left\langle \delta,\delta\right\rangle =-2$ can be
realized as calibrated cycle. Then the restriction of the automorphic function
$\exp\left(  \Phi_{3,19}\left(  \tau\right)  \right)  $ on the line
$\mathbb{R}\operatorname{Im}\omega_{Y}\subset\mathbb{R}^{k}+iV^{+}\subset
Pic(Y)\otimes\mathbb{C}$ is a periodic function such that the coefficients
$a_{n}$ in front of $\frac{\exp\left(  -int\right)  }{1-\exp(-int)}$ are
integer such that $a_{n}$ is equal to calibrated cycles $\delta$ such that
$vol(\delta)=n.$
\end{theorem}

\textbf{Proof: }We will prove the following Lemma:

\begin{lemma}
\label{calc1}The 2-cycle $\delta\in T(Y)=Pic(Y)^{\bot}\subset H^{2}%
(Y,\mathbb{Z})$ on the K3 surface $Y$ such that $\left\langle \delta
,\delta\right\rangle =-2,$ and $\left\langle \delta,\operatorname{Im}%
\omega_{Y}\right\rangle >0$ can be realized as calibrated cycle.
\end{lemma}

\textbf{Proof: }We know that $\operatorname{Im}\omega_{Y}\in H^{2,0}%
(Y,\mathbb{C})\oplus H^{0,2}(Y,\mathbb{C})\subset T(Y)\otimes\mathbb{R}.$ Let
us choose a CY metric $g$ on Y such that
\[
\left\langle \operatorname{Im}g,\delta\right\rangle =0\text{ and }\left\langle
\operatorname{Im}g,\omega_{Y}\right\rangle =0.
\]
Let us consider isometric deformation of $Y$ with respect to the CY metric
$g.$ From the properties of the isometric deformation of CY metrics on K3
surfaces studied in \cite{To80} we can change the complex structure on $Y$ in
such a way that \textbf{1. }$\left\langle \alpha\operatorname{Re}\omega
_{Y}+\beta\operatorname{Im}\omega_{Y},\delta\right\rangle >0$ for some real
numbers $\alpha$ and $\beta,$ \textbf{2. }the vector $\gamma\operatorname{Re}%
\omega_{Y}+\mu\operatorname{Im}\omega_{Y}$ in the three dimensional subspace
in $H^{2}(Y,\mathbb{R})$ spanned by $\operatorname{Re}\omega_{Y},$
$\operatorname{Im}\omega_{Y}$ and $\operatorname{Im}g$ perpendicular to
$\alpha\operatorname{Re}\omega_{Y}+\beta\operatorname{Im}\omega_{Y}\ $is such
that
\[
\left\langle \gamma\operatorname{Re}\omega_{Y}+\mu\operatorname{Im}\omega
_{Y},\delta\right\rangle =0
\]
and \textbf{3. }$\operatorname{Im}g$ and $\alpha\operatorname{Re}\omega
_{Y}+\beta\operatorname{Im}\omega_{Y}$ will be realized as the imaginary part
of a CY metric with respect to the new complex structure on $Y.$ It is easy to
see that the Poincare dual class of cohomology of $\delta$ can be realized as
a form of type $(1,1)$ with respect to the new isometric complex structure on
$Y.$ Thus as it was proved in \cite{PS} $\delta$ can be realized as a rational
non singular curve on the new K3 surface. Then the volume form of the
restriction of CY metric with imaginary part $\alpha\operatorname{Re}%
\omega_{Y}+\beta\operatorname{Im}\omega_{Y}$ on the rational curve with class
of homology $\delta$ will be $\operatorname{Im}\omega_{Y}.$ Lemma \ref{calc1}
is proved. $\blacksquare$

\begin{lemma}
\label{calc2}The restriction of the automorphic function $\exp\left(
\Phi_{3,19}\left(  \tau\right)  \right)  $ on the line $\mathbb{R}%
\operatorname{Im}\omega_{Y}\subset\mathbb{R}^{k}+iV^{+}\subset Pic(Y)\otimes
\mathbb{C}$ is a periodic function such that the coefficients $a_{n}$ in front
of $\frac{\exp\left(  -int\right)  }{1-\exp(-int)}$ are integer such that
$a_{n}$ is equal to calibrated cycles $\delta$ such that $vol(\delta)=n.$
\end{lemma}

\textbf{Proof:} The proof of Lemma \ref{calc2} is exactly the same as the
proof of Theorem \ref{to05}. $\blacksquare$ Theorem \ref{Bmodel} is proved.
$\blacksquare$

\begin{remark}
Theorem \ref{Bmodel} can be reformulated as follows: It is a well known fact
that the $\delta\in\Lambda_{K3}$ such that $\left\langle \delta,\delta
\right\rangle =-2$ can be realized as vanishing invariant cycles with a
monodromy group $r_{\delta}(v)=v+\left\langle v,\delta\right\rangle \delta$
for any $v\in\Lambda_{K3}.$ This means that there exists a family of K3
surfaces $\pi:\mathcal{X}\rightarrow D$ such that $X_{0}=\pi^{-1}(0)$ has an
isolated singularity of type $A_{n},D_{n}$ or $E_{6},$ $E_{7}$ and $E_{8}$ and
the monodromy operator acting on $H_{2}(X_{t},\mathbb{Z})$ by the reflection
$r_{\delta}$ described above. Thus in the B-model the partition function count
invariant calibrated vanishing cycle with a given volume when we choose the
complex structure such on the K3 surface with unimodular Picard group such
that $\operatorname{Im}\omega_{Y}\in H_{2}(X_{t},\mathbb{Z}).$
\end{remark}

\begin{conjecture}
The analogue of Theorem \ref{Bmodel} holds for the B-model of CY threefolds,
i.e. the partition function counts invariant vanishing calibrated cycles in
the B-model when the monodromy operator is of infinite order.
\end{conjecture}

\section{The Canonical Class of the Moduli Space of Polarized Algebraic K3}

\subsection{The Projection Formula}

\begin{theorem}
\label{prim}Let $l\in\Lambda_{K3}$ be a primitive vector such that
$\left\langle l,l\right\rangle =2n>0.$ Let $\left(  l\right)  ^{\perp}$ be the
sublattice in $\Lambda_{K3}$ orthogonal to $\mathbb{Z}l.$ Then we have$\left(
l\right)  ^{\perp}\approxeq\mathbb{Z}l^{\ast}\oplus\mathbb{U}^{2}\oplus\left(
-E_{8}\right)  ^{2},$ where $l^{\ast}$ is a primitive vector in $\Lambda_{K3}$
such that $\left\langle l^{\ast},l^{\ast}\right\rangle =-2n<0.$
\end{theorem}

\textbf{Proof: }According to \cite{PS} the subgroup $\mathcal{O}_{\Lambda
_{K3}}^{+}$ of index two that preserve the spinor norm acts transitively on
the primitive vectors with a fixed positive self intersection. Let us fix
$\mathbb{U}$ in $\Lambda_{K3}$ with a basis $e_{0}$ and $e_{1}$ such that
$\left\langle e_{i},e_{i}\right\rangle =0$ and $\left\langle e_{1}%
,e_{2}\right\rangle =1.$ Then $l=e_{1}+ne_{2}\in\mathbb{U}$ is a primitive
vector such that $\left\langle l,l\right\rangle =2n>0.$ Let $l^{\ast}%
=e_{1}-ne_{2}\in\mathbb{U}.$ Clearly $l^{\ast}$ is a primitive vector such
that $\left\langle l,l^{\ast}\right\rangle =0$ and $\left\langle l^{\ast
},l^{\ast}\right\rangle =-\left\langle l,l\right\rangle =-2n.$ Then we have
\begin{equation}
\left(  l\right)  ^{\perp}\approxeq\mathbb{Z}l^{\ast}\oplus\mathbb{U}%
\oplus\mathbb{U}\oplus\mathbb{E}_{8}(-1)\oplus\mathbb{E}_{8}(-1). \label{exp}%
\end{equation}
Theorem \ref{prim} is proved. $\blacksquare$

\begin{notation}
Let $\Lambda_{K3.n}:=\mathbb{Z}l^{\ast}\oplus\mathbb{U}^{2}\oplus\left(
-E_{8}\right)  ^{2}$ where $\left\langle l^{\ast},l^{\ast}\right\rangle =-2n.$
Let $\left\{  e_{1},e_{2},f_{1},f_{2},g_{1}\text{ and }g_{2}\right\}  $ be a
basis of $\mathbb{U}\oplus\mathbb{U}\oplus\mathbb{U}$ in $\left(
\ref{exp}\right)  $ such that $\left\langle e_{i},e_{i}\right\rangle
=\left\langle f_{i},f_{i}\right\rangle =\left\langle g_{i},g_{i}\right\rangle
=0$ and $\left\langle e_{1},e_{2}\right\rangle =\left\langle f_{1}%
,f_{2}\right\rangle =\left\langle g_{1},g_{2}\right\rangle =1.$
\end{notation}

\begin{theorem}
\label{Pr}The orthogonal projection of the discriminant $\mathfrak{D}_{3,19}$
on $\mathcal{M}_{K3,pa}^{2d}$ is $\mathfrak{D}_{n},$ where $\mathfrak{D}_{n}$
is the divisor in $\mathcal{M}_{K3,pa}^{2d}$ defined by the hyperplanes in
$\mathfrak{h}_{2,n}$ orhogonal to $l^{\ast}$ and all vectors $\delta\in
\Lambda_{n}$ such that $\left\langle \delta,\delta\right\rangle =-2.$
\end{theorem}

\textbf{Proof: }The proof of Theorem \ref{Pr} follows directly from the
following Lemma:

\begin{lemma}
\label{0} Let $\delta\in\Lambda_{K3}$ be such that $\left\langle \delta
,\delta\right\rangle =-2.$ Suppose that $\delta\notin\Lambda_{K3,n}.$ Then
there exists an automorphism $\sigma$ of the lattice $\Lambda_{K3.n}$ such
that $\Pr_{\mathbb{U}}\sigma(\delta)=l^{\ast}.$
\end{lemma}

\textbf{Proof: }The proof of Lemma \ref{0} is based on the following Propositions:

\begin{proposition}
\label{1} Let $\Lambda_{K3}=\mathbb{U\oplus L}$ and let $e_{1}$ and $e_{2}$ be
the isotropic generators of $\mathbb{U}.$ Let $l=e_{1}+ne_{2}\in\mathbb{U}$
and $n>0.$ Suppose that $\left\langle \delta,\delta\right\rangle =-2.$ Then
there exists an element $\sigma\in\Gamma_{n}$ such that in the representation
\[
\sigma(\delta)=n_{1}e_{1}+n_{2}e_{2}+\mu_{\sigma(\delta)},\mu_{\sigma(\delta)}%
\]
satisfies
\begin{equation}
\left\langle \Pr{}_{\mathbb{U}}(\sigma(\delta)),\Pr{}_{\mathbb{U}}%
(\sigma(\delta))\right\rangle <0\Longleftrightarrow\left\langle \mu
_{\sigma(\delta)},\mu_{\sigma(\delta)}\right\rangle >0. \label{exp0}%
\end{equation}

\end{proposition}

\textbf{Proof: }Let\textbf{ }$\delta=m_{1}e_{1}+m_{2}e_{2}+\mu_{\delta}.$ Let
us consider
\begin{equation}
\delta_{1}=k_{\delta_{1}}l^{\ast}+\mu_{\delta_{1}}\in\Lambda_{K3,n}%
,\left\langle \delta_{1},\delta_{1}\right\rangle =-2. \label{Exp0}%
\end{equation}
Then $\mu_{\delta_{1}}\in\left(  l^{\ast}\right)  ^{\perp}=\mathbb{L=U}%
\oplus\mathbb{U\oplus E}(-1)\mathbb{\oplus E}(-1).$ Since $\left\langle
l^{\ast},l^{\ast}\right\rangle <0,$ $\left(  \ref{Exp0}\right)  $ and
$\left\langle \delta_{1},\delta_{1}\right\rangle =-2$ $\ $then $\left\langle
\mu_{\delta_{1}},\mu_{\delta_{1}}\right\rangle >0.$ Let us consider the
reflection map
\[
\sigma(\delta)=r_{\delta_{1}}(\delta)=\delta+\left\langle \delta,\delta
_{1}\right\rangle \delta_{1},
\]
where $v\in\Lambda_{K3,n}.$ Let us compute the projection $\Pr{}_{\mathbb{U}%
}(\sigma(\delta))$ of $\sigma(\delta)$ on $\mathbb{U}$ spanned by $e_{1}$ and
$e_{2}.$ Direct computations show that%
\[
\Pr{}_{\mathbb{U}}(2n\delta)=2nm_{1}e_{1}+2nm_{2}e_{2}=nm_{1}\left(
l+l^{\ast}\right)  +m_{2}\left(  l-l^{\ast}\right)  =
\]%
\begin{equation}
\left(  nm_{1}+m_{2}\right)  l+\left(  nm_{1}-m_{2}\right)  l^{\ast}.
\label{Exp0a}%
\end{equation}
Direct computations show that%
\[
\Pr{}_{\mathbb{U}}(2n\sigma(\delta))=
\]%
\[
\left(  nm_{1}+m_{2}\right)  l-2n\left(  k_{\delta_{1}}\left\langle
\delta,\delta_{1}\right\rangle +\left(  nm_{1}-m_{2}\right)  \right)  l^{\ast
}.
\]
Suppose that $nm_{1}-m_{2}\neq0.$ So
\[
\left\langle \Pr{}_{\mathbb{U}}(2n\sigma(\delta)),\Pr{}_{\mathbb{U}}%
(2n\sigma(\delta))\right\rangle =
\]%
\[
\left(  nm_{1}+m_{2}\right)  ^{2}\left\langle l,l\right\rangle +\left(
2n\left(  k_{\delta_{1}}\left\langle \delta,\delta_{1}\right\rangle +\left(
nm_{1}-m_{2}\right)  \right)  \right)  ^{2}\left\langle l^{\ast},l^{\ast
}\right\rangle =
\]%
\begin{equation}
n\left(  nm_{1}+m_{2}\right)  ^{2}-n\left(  2n\left(  k_{\delta_{1}%
}\left\langle \delta,\delta_{1}\right\rangle +\left(  nm_{1}-m_{2}\right)
\right)  \right)  ^{2}. \label{expa}%
\end{equation}
We can choose $\delta_{1}$ such that $|k_{\delta_{1}}|$ is big enough. Thus
$\left(  \ref{expa}\right)  $ will imply $\left(  \ref{exp0}\right)  $%
\[
\left\langle \Pr{}_{\mathbb{U}}(2n\sigma(\delta)),\Pr{}_{\mathbb{U}}%
(2n\sigma(\delta))\right\rangle =
\]%
\[
n\left(  nm_{1}+m_{2}\right)  ^{2}-n\left(  2n\left(  k_{\delta_{1}%
}\left\langle \delta,\delta_{1}\right\rangle +\left(  nm_{1}-m_{2}\right)
\right)  \right)  ^{2}<0.
\]

Suppose that $nm_{1}-m_{2}=0.$ Then $\left(  \ref{Exp0a}\right)  $ implies
\[
\delta=\left(  m_{1}n+m_{2}\right)  l+\mu_{\delta}.
\]
Thus $\Pr{}_{\mathbb{U}}(\delta)=\left(  m_{1}n+m_{2}\right)  l.$ Let us
choose $\delta_{1}=k_{\delta_{1}}l^{\ast}+\mu_{\delta_{1}},$ such that
$\left\langle \delta_{1},\delta_{1}\right\rangle =-2,$ and $\left\langle
\delta,\delta_{1}\right\rangle \neq0.$ Let us compute
\begin{equation}
r_{\delta_{1}}(2n\delta)=2n\delta+2n\left\langle \delta,\delta_{1}%
\right\rangle \delta_{1}=\left(  m_{1}n+m_{2}\right)  l+\left\langle
\delta,\delta_{1}\right\rangle \left(  k_{\delta_{1}}l^{\ast}+\mu_{\delta_{1}%
}\right)  . \label{expb}%
\end{equation}
Thus $\left(  \ref{expb}\right)  $ implies that%
\[
\left\langle \Pr{}_{\mathbb{U}}(2nr_{\delta_{1}}(\delta)),\Pr{}_{\mathbb{U}%
}(2nr_{\delta_{1}}(\delta))\right\rangle =
\]%
\[
\left\langle \left(  m_{1}n+m_{2}\right)  l,\left(  m_{1}n+m_{2}\right)
l\right\rangle +\left\langle \left\langle \delta,\delta_{1}\right\rangle
k_{\delta_{1}}l^{\ast},\left\langle \delta,\delta_{1}\right\rangle
k_{\delta_{1}}l^{\ast}\right\rangle =
\]%
\[
\left(  m_{1}n+m_{2}\right)  ^{2}\left\langle l,l\right\rangle +\left(
\left\langle \delta,\delta_{1}\right\rangle k_{\delta_{1}}\right)
^{2}\left\langle l^{\ast},l^{\ast}\right\rangle =
\]%
\begin{equation}
\left(  m_{1}n+m_{2}\right)  ^{2}n-n\left(  \left\langle \delta,\delta
_{1}\right\rangle k_{\delta_{1}}\right)  ^{2}. \label{expba}%
\end{equation}
If we choose $\delta_{1}$ such that $\left\langle \delta,\delta_{1}%
\right\rangle \neq0$ and $\left\vert k_{\delta_{1}}\right\vert $ big is enough
then $\left(  \ref{expba}\right)  $ implies%

\[
\left\langle \Pr{}_{\mathbb{U}}(r_{\delta_{1}}(\delta)),\Pr{}_{\mathbb{U}%
}(2nr_{\delta_{1}}(\delta))\right\rangle =
\]
\[
\left(  m_{1}n+m_{2}\right)  ^{2}n-n\left(  \left\langle \delta,\delta
_{1}\right\rangle k_{\delta_{1}}\right)  ^{2}<0.
\]
Proposition \ref{1} is proved. $\blacksquare$

\begin{proposition}
\label{min}Let $\Gamma_{n}$ be the generated by the reflections
\[
r_{\delta}:v\rightarrow\left\langle v,\delta\right\rangle \delta
\]
for
\[
\delta\in\mathbb{Z}l^{\ast}\oplus\mathbb{U}\oplus\mathbb{U}\oplus
\mathbb{E}_{8}(-1)\oplus\mathbb{E}_{8}(-1).
\]
Let $\delta_{\min}\in\left\{  \Gamma_{n}\delta\right\}  $ be such that
\begin{equation}
\left\langle \mu_{\delta_{\min}},\mu_{\delta_{\min}}\right\rangle
=\underset{\left\{  \sigma\in G_{n}\right\}  }{\min}\left\langle \mu
_{\sigma(\delta)},\mu_{\sigma(\delta)}\right\rangle \geq0. \label{Min}%
\end{equation}
Then $\left\langle \mu_{\delta_{\min}},\mu_{\delta_{\min}}\right\rangle =0.$
\end{proposition}

\textbf{Proof: }Let $\delta_{\min}=pe_{1}+qe_{2}+\mu_{\delta_{\min}}.$ Then
acoording to $\left(  \ref{exp0a}\right)  $ we have
\begin{equation}
\Pr{}_{\mathbb{U}}\left(  2n\delta_{\min})\right)  =\left(  pn+q\right)
l+\left(  pn-q\right)  l^{\ast}, \label{Expob}%
\end{equation}
where $l=e_{1}+ne_{2}$ and $l^{\ast}=e_{1}-ne_{2}.$ So $\left(  \ref{Min}%
\right)  $ implies
\[
\left\langle \Pr{}_{\mathbb{U}}\left(  2n\delta_{\min})\right)  ,\Pr
{}_{\mathbb{U}}\left(  2n\delta_{\min})\right)  \right\rangle <0
\]
which is equivalent to
\begin{equation}
\left(  pn+q\right)  ^{2}-\left(  pn-q\right)  ^{2}<0. \label{Expod}%
\end{equation}
Let us choose
\begin{equation}
\left\langle \delta,\delta\right\rangle =-2\text{ and }\delta=k_{\delta
}l^{\ast}+\mu_{\delta}. \label{-2}%
\end{equation}
Let us consider $r_{\delta}(\delta_{\min})=\delta_{\min}+\left\langle
\delta,\delta_{\min}\right\rangle \delta.$ Direct computaions using $\left(
\ref{-2}\right)  $ and $\left(  \ref{expa}\right)  $ show that
\begin{equation}
r_{\delta}(2n\delta_{\min})=\left(  pn+q\right)  l+\left(  \left(
pn-q\right)  +2nk_{\delta}\left\langle \delta_{\min},\delta\right\rangle
\right)  l^{\ast}+\mu_{r_{\delta_{1}}(\delta)}. \label{-1}%
\end{equation}

\begin{remark}
\label{C}Let $\delta=k_{\delta}l^{\ast}+\mu_{\delta}$ and $\delta
_{1}=-k_{\delta}l^{\ast}+\mu_{\delta_{1}}$satisfy $\left\langle \delta
,\delta\right\rangle =\left\langle \delta_{1},\delta_{1}\right\rangle =-2$
then we can choose $\mu_{\delta_{1}}$ to be such that the sign of
$\left\langle \delta_{\min},\delta\right\rangle $ to be the same as that of
$\left\langle \delta_{\min},\delta_{1}\right\rangle .$
\end{remark}

\textbf{Proof}: Let $f_{i}$ and $g_{i}$ be the generators of $\mathbb{U}%
\oplus\mathbb{U},$ where $\left\langle f_{i},f_{i}\right\rangle =\left\langle
g_{i},g_{i}\right\rangle =0$ and $\left\langle f_{1},f_{2}\right\rangle
=\left\langle g_{1},g_{2}\right\rangle =1.$ Then we can choose
\begin{equation}
\mu_{\delta_{\min}}=f_{1}+\frac{\left\langle \mu_{\delta_{\min}},\mu
_{\delta_{\min}}\right\rangle }{2}f_{2},\text{ }\mu_{\delta_{1}}=g_{1}%
+\frac{\left\langle \mu_{\delta_{1}},\mu_{\delta_{1}}\right\rangle }{2}%
g_{2}+mf_{2}. \label{A}%
\end{equation}
Then it is clear that
\[
\left\langle \delta_{\min},\delta\right\rangle =\left\langle pe_{1}%
+qe_{2},k_{\delta}\left(  e_{1}+ne_{2}\right)  \right\rangle +\left\langle
\mu_{\delta_{\min}},\mu_{\delta}\right\rangle =
\]%
\[
k_{\delta}q+\left\langle \mu_{\delta_{\min}},\mu_{\delta}\right\rangle .
\]
On the other hand we derive from $\left(  \ref{A}\right)  $
\begin{equation}
\left\langle \delta_{\min},\delta_{1}\right\rangle =-k_{\delta}q+\left\langle
\mu_{\delta_{\min}},\mu_{\delta}\right\rangle =-k_{\delta}q+m. \label{B}%
\end{equation}
It is clear that we can choose $m$ such that the sign of $\left\langle
\delta_{\min},\delta\right\rangle $ to be the same as the sign of
$\left\langle \delta_{\min},\delta_{1}\right\rangle .$ Remark \ref{C} is
proved. $\blacksquare$

Thus Remark \ref{C} implies that we can choose the sign of $k_{\delta}$ in the
expression of $\delta$ such that the sign of $k_{\delta}\left\langle
\delta_{\min},\delta\right\rangle $ to be the oposite of the sign of $\left(
pn-q\right)  .$ So
\begin{equation}
\left(  \left(  pn-q\right)  +2nk_{\delta}\left\langle \delta_{\min}%
,\delta\right\rangle \right)  ^{2}<\left(  pn-q\right)  ^{2}. \label{min3c}%
\end{equation}
Thus $\left(  \ref{Expod}\right)  $ and $\left(  \ref{min3c}\right)  $ imply
that%
\[
\left\langle \Pr{}_{\mathbb{U}}\left(  2nr_{\delta}(\delta_{\min})\right)
,\Pr{}_{\mathbb{U}}\left(  2nr_{\delta}(\delta_{\min})\right)  \right\rangle
=
\]%
\begin{equation}
4n^{2}\left(  \left(  pn+q\right)  ^{2}-\left(  \left(  pn-q\right)
+2nk_{\delta}\left\langle \delta_{\min},\delta\right\rangle \right)
^{2}\right)  <0. \label{min3a}%
\end{equation}
So $\left(  \ref{min3a}\right)  $ implies that
\begin{equation}
\left\langle \mu_{r_{\delta}}(\delta_{\min}),\mu_{r_{\delta}}(\delta_{\min
})\right\rangle \geq0. \label{-0}%
\end{equation}
Since
\[
2nr_{\delta}(\delta_{\min})=\Pr{}_{\mathbb{U}}\left(  2nr_{\delta}%
(\delta_{\min})\right)  +\mu_{r_{\delta_{1}}(\delta)}%
\]
then $\left(  \ref{-0}\right)  ,$ $\left(  \ref{min3c}\right)  $ and $\left(
\ref{min3a}\right)  $ show that%
\[
\left\langle 2n\mu_{\min},2n\mu_{\min}\right\rangle =-8n^{2}+4n^{2}\left(
\left(  pn-q\right)  ^{2}-\left(  pn+q\right)  ^{2}\right)  >
\]%
\[
-8n^{2}+4n^{2}\left(  \left(  \left(  pn-q\right)  +2nk_{\delta}\left\langle
\delta_{\min},\delta\right\rangle \right)  ^{2}-\left(  pn+q\right)
^{2}\right)  =
\]%
\[
\left\langle 2n\mu_{r_{\delta}(\delta_{\min})},2n\mu_{r_{\delta}(\delta_{\min
})}\right\rangle >0.
\]
So we get that%
\begin{equation}
\left\langle \mu_{\min},\mu_{\min}\right\rangle >\left\langle \mu_{r_{\delta
}(\delta_{\min})},\mu_{r_{\delta}(\delta_{\min})}\right\rangle \geq0.
\label{Min1a}%
\end{equation}
Thus we get a contradiction with $\left\Vert \mu_{\delta_{\min}}\right\Vert
^{2}>0$ being the minimal value. Proposition \ref{min} is proved.
$\blacksquare$ Proposition \ref{min} implies Lemma \ref{0}. $\blacksquare$
Lemma \ref{0} implies Theorem \ref{Pr}. $\blacksquare$

\subsection{The Divisor of the Restriction of the Automorphic Form on
$\mathfrak{M}_{K3,n}$}

Let us consider the moduli space $\mathfrak{M}_{K3,n}$ of pseudo polarized
algebraic K3 surfaces with a polarization class $l\in\Lambda_{K3},$ where $l$
is a primitive vector in $\Lambda_{K3}$ such that $\left\langle
l,l\right\rangle =2n.$ Then according to \cite{PS} and \cite{D} we have
$\mathfrak{M}_{K3,n}=\Gamma_{n}\backslash\mathfrak{h}_{2,19},$ where
$\Gamma_{n}:=\left\{  \phi\in\mathcal{O}_{\Lambda_{K3}}^{+}|\phi(l)=l\right\}
.$ According to \cite{To80} we can define $\mathfrak{h}_{2,19}$ as one of the
open components of the quadric $\mathcal{Q}\subset\mathbb{P}(\Lambda
_{K3,n}\otimes\mathbb{C})$ defined as follows%
\[
\mathcal{Q}:=\left\{  u\in\mathbb{P}(\Lambda_{K3,n}\otimes\mathbb{C}%
)|\left\langle u,u\right\rangle =0\text{ and }\left\langle u,\overline
{u}\right\rangle >0.\right\}
\]
Let us define $\mathfrak{D}_{n}$ in $\mathfrak{M}_{K3,n}$ as follows: Let
$\lambda\in\Lambda_{K3,n},$ then%
\[
\mathcal{H}_{\lambda}:=\left\{  u\in\mathbb{P}(\Lambda_{K3,n}\otimes
\mathbb{C})|\left\langle u,\lambda\right\rangle =0.\right\}
\]
Let
\begin{equation}
\mathcal{D}_{n}:=\left(
{\displaystyle\bigcup\limits_{\left\langle \delta,\delta\right\rangle
=-2\text{ \& }\delta\in\Lambda_{K3,n}}}
\left(  \mathfrak{h}_{2,19}\cap\mathcal{H}_{\delta}\right)  \right)
\cup\left(
{\displaystyle\bigcup\limits_{\phi\in\Gamma_{n}}}
\left(  \mathfrak{h}_{2,19}\cap\mathcal{H}_{\phi(l^{\ast})}\right)  \right)
.\label{deco}%
\end{equation}
Then $\mathfrak{D}_{n}:=\Gamma_{n}\backslash\mathcal{D}_{n}.$

\begin{theorem}
\label{pc}There exists an automorphic form $\Psi_{19,n}$ on $\mathfrak{M}%
_{K3,n}=\Gamma_{n}\backslash\mathfrak{h}_{2,19}$ such that the zero set of
$\Psi_{19,n}$ is $\mathfrak{D}_{n}.$
\end{theorem}

\textbf{Proof: }According to the results of Harvey, Moore and Borcherds on we
can find an automorphic form $\left\vert \Psi_{\Lambda_{K3}}\right\vert ^{2}$
on the moduli space of Einstein metrics $\mathcal{O}_{\Lambda_{K3}}^{+}$%
$\backslash$%
$\mathfrak{h}_{3,19}$ such that its zeros are exactly on the discriminant
locus of $\mathcal{O}_{\Lambda_{K3}}^{+}$%
$\backslash$%
$\mathfrak{h}_{3,19}.$ recall that the discriminant locus on $\mathcal{O}%
_{\Lambda_{K3}}^{+}$%
$\backslash$%
$\mathfrak{h}_{3,19}$ is defined as the set of three dimensional positive
vector subspaces in $\Lambda_{K3}\otimes\mathbb{R}$ perpendicular to $\delta$
such that $\left\langle \delta,\delta\right\rangle =-2$ modulo the action of
the arithmetic group $\mathcal{O}_{\Lambda_{K3}}^{+}$. The moduli space
$\mathfrak{M}_{K3,n}=\Gamma_{n}\backslash\mathfrak{h}_{2,19}$ can be embedded
in $\mathcal{O}_{\Lambda_{K3}}^{+}$%
$\backslash$%
$\mathfrak{h}_{3,19}$ as the set of all three dimensional oriented subspaces
in $\Lambda_{K3}\otimes\mathbb{R}$ containing the polarization vector $l$
modulo the action of $\mathcal{O}_{\Lambda_{K3}}^{+}.$ The restriction of some
power of $\Psi_{\Lambda_{K3}}$ on $\mathfrak{M}_{K3,n}$ will give us an
automorphic form $\Psi_{19,n}$ on $\mathfrak{M}_{K3,n}.$ Thus we have the
following obvious fact:

\begin{remark}
\label{pc1} The zero set of the restriction of $\Psi_{\Lambda_{K3}}%
=\exp\left(  \Phi_{\Lambda_{3,19}}(\tau)\right)  $ on $\mathfrak{M}_{K3,n}$ is
the projection of the zero set of $\Psi_{\Lambda_{K3}}=\exp\left(
\Phi_{\Lambda_{3,19}}(\tau)\right)  $ on $\mathfrak{M}_{K3,n}.$
\end{remark}

Thus we need to compute the projection of the zero set of $\exp\left(
\Phi_{\Lambda_{3,19}}\right)  $ on $\Gamma^{+}\backslash\mathfrak{h}_{3,19}$
to $\mathfrak{M}_{K3,n}=\Gamma_{n}\backslash\mathfrak{h}_{2,19}.$ Theorem
\ref{pc} will follow from the following Lemma:

\begin{lemma}
\label{3} The zero set of $\Psi_{19,n}$ on $\mathfrak{M}_{K3,n}$ is
$\mathfrak{D}_{n}.$
\end{lemma}

\textbf{Proof: }Let $\delta\in\Lambda_{K3}$ be such that $\left\langle
\delta,\delta\right\rangle =-2.$ Let $\Pr_{l,n}(\delta)\in\Lambda_{K3,n}$ be
the orthogonal projection of $\delta$ on $\Lambda_{K3,n}.$ If
\[
\Pr{}_{l,n}(\delta)=\delta\iff\left\langle l,\delta\right\rangle =0,
\]
then it implies that the component $%
{\displaystyle\bigcup\limits_{\left\langle \delta,\delta\right\rangle
=-2\text{ \& }\delta\in\Lambda_{K3,n}}}
\left(  \mathfrak{h}_{2,19}\cap\mathcal{H}_{\delta}\right)  $ in the
expression $\left(  \ref{deco}\right)  $ defines the components of
$\mathfrak{D}_{n}:=\Gamma_{n}\backslash\mathcal{D}_{n}$ corresponding to the
vectors with $-2$ norm in $\Lambda_{K3,n}.$

Suppose that $\delta\in\Lambda_{K3},$ $\left\langle \delta,\delta\right\rangle
=-2$ and $\Pr{}_{l,n}(\delta)\neq\delta.$ Theorem \ref{Pr} implies that we can
find $\sigma\in\Gamma_{n}$ such that $\sigma(\delta)=m_{1}e_{1}+m_{2}e_{2}.$
Thus $\Pr{}_{l,n}(\delta)=k_{\delta}l^{\ast}.$ Then
\begin{equation}
\pi\left(  \mathcal{H}_{\delta}\cap\mathfrak{h}_{2,19}\right)  =\pi\left(
H_{l^{\ast}}\cap\mathfrak{h}_{2,19}\right)  \label{min5}%
\end{equation}
where $\pi:\mathfrak{h}_{2,19}\rightarrow\Gamma_{n}\backslash\mathfrak{h}%
_{2,19}.$Thus $\left(  \ref{min5}\right)  $ implies Lemma \ref{3}.
$\blacksquare$

Theorem \ref{pc} is proved. $\blacksquare$

\begin{corollary}
\label{JT}The zero set of the restriction of the automorphic form
$\Psi_{\Lambda_{K3}}=\exp\left(  \Phi_{\Lambda_{3,19}}(\tau)\right)  $ on
$\mathfrak{M}_{K3,n}$ is a divisor $\mathfrak{D}_{n}:=\Gamma_{n}%
\backslash\mathcal{D}_{n}$ which consists of two components $\pi\left(
H_{l^{\ast}}\cap\mathfrak{h}_{2,19}\right)  $ and $\pi\left(  H_{\delta}%
\cap\mathfrak{h}_{2,19}\right)  ,$ where $\delta\in\left(  l^{\ast}\right)
^{\bot}=\mathbb{U}\oplus\mathbb{U}\oplus\mathbb{E}_{8}(-1)\oplus\mathbb{E}%
_{8}(-1)$ and $\pi:\mathfrak{h}_{2,19}\rightarrow\Gamma_{n}\backslash
\mathfrak{h}_{2,19}=\mathfrak{M}_{K3,n}.$
\end{corollary}

\textbf{Proof: }Corollary \ref{JT} follows from Theorem \ref{Bor} which
implies that the divisor $\pi\left(  H_{\delta}\cap\mathfrak{h}_{2,19}\right)
$ is an irreducible since we assumed that $\delta\in\left(  l^{\ast}\right)
^{\bot}=\mathbb{U}\oplus\mathbb{U}\oplus\mathbb{E}_{8}(-1)\oplus\mathbb{E}%
_{8}(-1)$. The irreduciblity of $\pi\left(  H_{l^{\ast}}\cap\mathfrak{h}%
_{2,19}\right)  $ follows from Theorem  \ref{Pr}. $\blacksquare$

Corollary \ref{JT} is generalization of the results obtained in \cite{JT},
\cite{JT0} and \cite{JT1}.

\end{document}